\documentclass[11pt,letter,reqno]{amsart}

\usepackage{amstext,amsmath,amssymb,amsthm}
\usepackage[dvips]{graphics,color}
\usepackage{epsfig}

\usepackage[pdftex]{geometry}
\usepackage{graphicx}
\numberwithin{equation}{section}

\newtheorem{thm}{Theorem}[section]
\newtheorem{lma}[thm]{Lemma}

\newtheorem{rem}[thm]{Remark}

\newcommand{\p}{{\partial}}

\newcommand{\Div}{{\rm div }}

\newcommand{\bld}[1]{\boldsymbol{#1}}

\newcommand{\bq}{\bld{q}}

\newcommand{\lbk}{\lbrack \! \lbrack}
\newcommand{\rbk}{\rbrack \! \rbrack}
\newcommand{\lb}{\lbrace \!\! \lbrace}
\newcommand{\rb}{\rbrace \!\! \rbrace}

\newcommand\ww[1]{\textcolor{black}{#1}}
\newcommand\xs[1]{\textcolor{black}{#1}}

\usepackage{tikz}
\usetikzlibrary{shapes.geometric, arrows.meta, positioning}

\tikzstyle{startstop} = [rectangle, rounded corners, 
minimum width=2cm, 
minimum height=1cm,
text centered, 
draw=black, 
fill=red!30]

\tikzstyle{io} = [rectangle, rounded corners,
trapezium stretches=true, 
trapezium left angle=70, 
trapezium right angle=110, 
minimum width=3cm, 
minimum height=1cm, text centered, 
draw=black, fill=blue!30]

\tikzstyle{process} = [rectangle, 
minimum width=3cm, 
minimum height=1cm, 
text centered, 
text width=5cm, 
draw=black, 
fill=orange!30]

\tikzstyle{process2} = [rectangle, 
minimum width=3cm, 
minimum height=1cm, 
text centered, 
text width=14cm, 
draw=black, 
fill=orange!30]

\tikzstyle{process3} = [rectangle, 
minimum width=3cm, 
minimum height=1cm, 
text centered, 
text width=12cm, 
draw=black, 
fill=orange!30]

\tikzstyle{decision} = [trapezium, 
minimum width=3cm, 
minimum height=1cm, 
text width=6cm, 
text centered, 
draw=black, 
fill=green!30]
\tikzstyle{decision2} = [trapezium, 
minimum width=3cm, 
minimum height=1cm, 
text width=8.5cm, 
text centered, 
draw=black, 
fill=green!30]
\tikzstyle{arrow} = [thick,->,>=stealth]

\begin{document}

\title[A high order discontinuous Galerkin method for EIT]{A high order discontinuous Galerkin method for the recovery of the conductivity in Electrical Impedance Tomography}




\author[Li \and Wang]{Xiaosheng Li\and Wei Wang}
\address{Department of Mathematics and Statistics,
Florida International University, Miami, FL 33199, U.S.A.}
\email{xli@fiu.edu}
\address{Department of Mathematics and Statistics,
Florida International University, Miami, FL 33199, U.S.A.}
\email{weiwang1@fiu.edu}

\maketitle

\begin{abstract}

In this work, we develop an efficient high order discontinuous Galerkin (DG) method for solving the Electrical Impedance Tomography (EIT). EIT is a highly nonlinear ill-posed inverse problem where the interior conductivity of an object is recovered from the surface measurements of voltage and current flux. 
 We first propose a new optimization  problem based on the recovery of the conductivity from the Dirichlet-to-Neumann map to minimize  the mismatch between the predicted current  and the measured current on the boundary. And we further prove the existence of the minimizer. 
 Numerically the optimization problem is solved by a third order DG method with quadratic polynomials.
 Numerical results for several two-dimensional problems with both single and multiple inclusions are demonstrated to show the high {accuracy and efficiency} of the proposed high order DG method.  Analysis and computation for discontinuous conductivities are also studied in this work. 
\vspace{.5\baselineskip}


\vspace{.5\baselineskip}


$\qquad\quad\ \ $Dirichlet-to-Neumann map

\end{abstract}





\section{Introduction}
Electrical Impedance Tomography (EIT) is an imaging method to find the conductivity
of an object by making current and voltage measurements at the boundary.
It has many applications including the early diagnosis
of breast cancer \cite{CIN,ZG}, detection
of pneumothorax \cite{CCGB}, monitoring pulmonary
functions \cite{INGC},  detection of leaks from buried pipes \cite{JGP} and in underground storage tanks \cite{RDBLR}, as well as many industrial applications \cite{W}.
EIT is a typical inverse boundary value problem. The unique determination results have been obtained in \cite{AP,Na,SU}. The stability estimates \cite{Al,AV,BFR,Ma} indicate that such inverse problem is
severely ill-posed. We refer to Uhlmann's survey article \cite{U} for
the detailed development of the inverse boundary value problems in the theoretical aspect
since Calder\'{o}n's fundamental work \cite{C}.

Computationally, due to the high degree of nonlinearity and severe ill-posedness of the image problem, many efforts have been made in the development of 
efficient and stable numerical reconstruction algorithms.
These algorithms include the direct methods \cite{CIZ,KG,KLMS,SMI,Sy}, iterative methods \cite{CINSG,CCT,Do,Ha,Jin12,JM,KLZ,LR1}, variational methods \cite{BGZ,KV}, statistical approaches \cite{KKSV,KS}, neutral networks \cite{AG, BYZZ, FY}, among others. We refer to the survey articles \cite{Bo2,JMrev,KL,eit_error}. In practice, the full knowledge of the boundary measurement is not known. 
\xs{Only  the
data from a finite number of experiments is available, and the data may also contain some noises.}
The inverse problem is usually translated to an optimization problem to minimize the mismatch between the model predicted data and the measured data on the boundary. Because of the ill-posedness arising in the EIT problem, some types of regularization techniques \cite{EHN} are needed to stabilize the problem. 
\xs{The Tikhonov regularization method is widely recognized as the most commonly employed technique. 
The optimization problem can be solved iteratively, where inside each iteration the forward problem needs to be solved numerically. 
As the accuracy of the algorithm highly relies on the accuracy of the forward problem, an efficient and accurate forward solver is in desire.
There are many numerical techniques present to solve the forward problem, which can be
modeled by elliptic types of problems.}
Since finite volume and finite difference approaches  generally need regular grids,
the finite element method is commonly used for EIT applications.
\ww{Recent work using finite element method including discontinuous Galerkin,  stochastic Galerkin and weak Galerkin methods as a forward solver in the simulations for EIT includes  \cite{eit_forward2,Jin14, dgeit, 6, Zou20,Zou17,eit_hp,WG,eit_forward1}, etc.
However, there are not many works of high order methods to simulate both the forward and inverse of EIT problems.}


In this paper, we develop a high order discontinuous Galerkin (DG) method as the reconstruction method to solve the forward elliptic problems. 
The DG method is a class of finite element methods using completely discontinuous
piecewise polynomial space for the numerical solution and the test functions. 
An introduction of 
the development of DG methods can be found in the survey papers and books  \cite{dgbook, CSrev, CSp,Dp,
Hesthaven, Riviere}.
Recent developments, mainly for elliptic problems, include \cite{dg2, 1, dg3, CoDo07, hdg, hdg2,  wdg, Wang11, YS2}.
There
are several distinctive features that make DG attractive in applications, which include
the local conservativity, the ability for easily handling irregular meshes with hanging
nodes and boundary conditions, the flexibility for hp-adaptivity.
Besides those, 
DG also has advantages to deal with rough coefficients, especially the coefficients
containing discontinuities or multiscales.
And thus, DG methods have been well developed in
a wide range of applications.
However,  to the authors' best knowledge, there is little work for DG method in solving EIT problems.
In particular, it is difficult for traditional finite element methods to go high order in multidimensions because it requires continuities on the element boundaries.
Furthermore, it is also challenging for traditional finite element methods to deal with discontinuities such as in the conductivity coefficients.
\xs{Those advantages, the hp-adaptivity to go high order and the ability to deal with rough coefficients, make DG method attractive and suitable for EIT problems.} 
Thus, we would like to design a high order DG method and apply it to solve EIT problems.


In our work, we focus on the recovery of the conductivity from the Dirichlet-to-Neumann map, where the given voltage is applied on the boundary and the corresponding current flux after the interaction of the electromagnetic wave with the object is measured. We construct an optimization problem to minimize the mismatch between the predicted current from the Dirichlet-to-Neumann map and the measured current on the boundary with Tikhonov regularization. We prove the existence of the minimizer and derive the derivative formulas associated with the Dirichlet-to-Neumann map. We then apply our newly designed high order DG method to solve this EIT problem. 

This paper is organized as follows. In Section~2, we study the minimization problem for general conductivities, state the iteration procedure, and derive the formulas for the derivatives of the associated operators. In Section~3, we introduce the DG method for the forward problem. In Section~4, we describe the detailed algorithm for solving the inverse problem. Several numerical examples are presented to demonstrate the performance of the proposed method in Section~5. 
\ww{A special case of piecewise continuous conductivity is discussed in Section~6.}
In Section~7, we draw conclusions and make suggestions for further work.

\section{The minimization problem}

In this section we state the mathematical model and formulate the minimization problem. Suppose that $\Omega$ is a bounded and simply connected 
domain in $\mathbb{R}^d$ ($d\geq 2$) with Lipschitz boundary, and let the voltage potential $u$ 
solve the Dirichlet problem for the conductivity equation
\begin{equation}
\left\{\begin{array}{ll}\Div (\sigma\nabla u)=0&\quad \mbox{in}\
\Omega\\
u=f&\quad \mbox{on}\
\partial\Omega\end{array} \right.
\label{eq_int_u}
\end{equation}
where the conductivity function $\sigma$ is positive and bounded in $\Omega$. This problem has a unique solution $u\in H^1(\Omega)$ for any $f\in H^{\frac{1}{2}}(\partial\Omega)$ by the \xs{Lax-Milgram}
theorem. On the boundary, we can measure the outgoing current flux for a given boundary
voltage. The Dirichlet-to-Neumann map 
\begin{equation*}
F(\sigma, \cdot): H^{\frac{1}{2}}(\partial\Omega)\to H^{-\frac{1}{2}}(\partial\Omega)
\end{equation*}
is given by 
\begin{equation}
F(\sigma, f)=\left(\sigma\frac{\p u}{\p \nu}\right)\Big|_{\partial\Omega} \label{eq:f}
\end{equation}
where $\nu$ is the unit outer normal of $\partial\Omega$. The inverse problem consists of recovering $\sigma$ from $F(\sigma, \cdot)$.
We suppose that the conductivity is known on the boundary, and our main aim is to reconstruct the conductivity inside the domain. 

For the conductivity equation \eqref{eq_int_u}, when $f\in H^{\frac{1}{2}}(\partial\Omega)$, we know 
$F(\sigma, f)\in H^{-\frac{1}{2}}(\partial\Omega)$. It is inconvenient to compute with $H^{-\frac{1}{2}}(\partial\Omega)$ norm. In order to work with the $L^2(\partial\Omega)$ norm for easy computation, 
we need more regularity for the conductivity, the boundary data, and the domain, so that the regularity theory for elliptic equation can be used. This is new and different from EIT problem with Neumann-to-Dirichlet map \xs{(see Remark \ref{rem:n2d})}. 
Denote  
\begin{equation*}
\mathcal{A}=\{\sigma\in W^{1,\infty}(\Omega): 0<c_1<\sigma<c_2,\ |\nabla \sigma|< N, \mbox{ and } \sigma \mbox{ is known on } \partial\Omega\}\, 
\end{equation*}
the  admissible set for the conductivity, where $c_1$, $c_2$ and $N$ are fixed numbers. We suppose that $\Omega$ has $C^{1,1}$ boundary or 
$\Omega$ is a convex domain. If we take $f\in H^{\frac{3}{2}}(\partial\Omega)$, from elliptic theorem, then 
$F(\sigma, f)\in H^{\frac{1}{2}}(\partial\Omega)$.  We also endow 
$\mathcal{A}$ with the $H^1(\Omega)$ norm. 
\begin{rem}
When the conductivity is a piecewise continuous function, we can release the higher regularity requirement. This case is studied in Section~6.
\end{rem}
\begin{rem}\label{rem:n2d}
The Neumann-to-Dirichlet map 
$
G(\sigma, \cdot): H^{-\frac{1}{2}}(\partial\Omega)\to H^{\frac{1}{2}}(\partial\Omega)
$
is given by 
$
G(\sigma, g)=u|_{\partial\Omega}
$,
where $u$ is the solution of
\begin{equation*}
\left\{\begin{array}{ll}\Div (\sigma\nabla u)=0&\quad \mbox{in}\
\Omega\\
\sigma\frac{\p u}{\p \nu}=g&\quad \mbox{on}\
\partial\Omega\end{array} \right.
\end{equation*}
with $\int_{\partial\Omega}gds=1$. The $L^2(\partial\Omega)$ norm can be employed since $G(\sigma, g)\in H^{\frac{1}{2}}(\partial\Omega)\subset L^2(\partial\Omega)$.
\end{rem}

\xs{The Dirichlet-to-Neumann map involves an infinite number of boundary measurements. However, in practical applications, it is only feasible to collect a finite number of measurements, which may also contain noise. As a consequence of these measurement limitations, we can only obtain an approximate conductivity that deviates from the true conductivity. The accuracy of this approximation is contingent on the degree of noise present in the measurements.
Let $\sigma^{true}$ be the true conductivity we plan to reconstruct. Denote $f_j$ the imposed voltage on the boundary,  for $j=1, \cdots, M$ with $M$ being the number of experiments. Let $g_j^{true}=F(\sigma^{true},f_j)$ be the exact current flux on the boundary and $g_j^\delta$ be the measured current flux on the boundary, which contains some noises.}
So we have $M$ pairs of the available data $(f_j, g_j^\delta)$. The inverse problem we consider is to minimize the functional 
\begin{equation}
R(\sigma)=\frac{1}{2}\sum_{j=1}^{M}\|F(\sigma,f_j)-g_j^\delta\|_{L^2(\partial \Omega)}^2+\frac{\alpha}{2}\|\sigma-\sigma^0\|_{H^1(\Omega)}^2
\label{eq_int_R}
\end{equation}
over the admissible set $\mathcal{A}$. 
The first item describes the mismatch between model predictions and measurements. The second term is the regularization term, where $\alpha>0$ is the regularization parameter and $\sigma^0$ is the initial guess of the true conductivity. The minimizer is considered as an approximation to the true conductivity. 

\subsection{Existence of the minimizer.} We show that there exists at least one minimizer to the functional $R(\sigma)$. The proof is based on the continuity of $F(\sigma, f)$ for $\sigma\in\mathcal{A}$. The current literature is mainly for the Neumann-to-Dirichlet map \xs{(see, for example, \cite{CINSG,CCT,Do,Jin12,JM})}. Here we consider the Dirichlet-to-Neumann map. We need some regularity results for the solution to \eqref{eq_int_u}. From the standard elliptic theory, we first know that for $\sigma\in\mathcal{A}$ and $f\in H^{\frac{3}{2}}(\partial\Omega)$, we have $u\in H^2(\Omega)$ and
\begin{equation}
\|u\|_{H^2(\Omega)}\leq C \|f\|_{H^{\frac{3}{2}}(\partial\Omega)}\label{est_u_h2}
\end{equation}
where $C$ may depend on $c_1$, $c_2$, and $N$ in $\mathcal{A}$, but is independent of $u$ and $f$. Here and below we use $C$ to denote such generic constants, and they may vary from line to line. We also need the following results of Meyers's reverse H\"{o}lder estimates \cite{Me,GM,RS}. This result is also used in \cite{JM}.
\begin{thm} Suppose that $0<c_1<\sigma<c_2$ in $\Omega\in\mathbb{R}^d$ ($d\geq 2$). Let $u\in H^1(\Omega)$ be a weak solution of
\begin{equation*}
\Div (\sigma\nabla u)=\Div\, \mathbf{G} +h\quad \mbox{in}\ \ \Omega.
\end{equation*}
Then there exists $p>2$, depending on $c_1$, $c_2$ and $d$, such that $u\in W^{1,p}(\Omega)$ and
\begin{equation*}
\|u\|_{W^{1,p}(\Omega)}\leq C(\|u\|_{H^1(\Omega)}+\|\mathbf{G}\|_{L^p(\Omega)}+\|h\|_{L^p(\Omega)})
\end{equation*}
where $C$ depends on $c_1$, $c_2$, $\Omega$ and $p$. 
\label{thm_Meyers}
\end{thm}
Applying Theorem~\ref{thm_Meyers} to the solution of \eqref{eq_int_u}, we know that $u\in W^{1,p}(\Omega)$ for some $p>2$ and
\begin{equation}
\|u\|_{W^{1,p}(\Omega)}\leq C\|u\|_{H^1(\Omega)}\leq C\|f\|_{H^{\frac{1}{2}}(\partial\Omega)}\leq C\|f\|_{H^{\frac{3}{2}}(\partial\Omega)}.\label{est_u_w1p}
\end{equation}
Next we show $u\in W^{2,p}(\Omega)$. Denote $\mathbf{x}=(x_1,x_2,\cdots,x_d)$ and $w=\frac{\partial u}{\partial x_i}$ for some $1\leq i\leq d$. From $u\in H^2(\Omega)$ we know $w\in H^1(\Omega)$ and $w$ satisfies
\begin{equation*}
\Div (\sigma\nabla w)=-\Div (\frac{\partial \sigma}{\partial x_i}\nabla u)\quad \mbox{in}\ \ \Omega.
\end{equation*}
Applying Theorem~\ref{thm_Meyers} to the above equation, we obtain $w\in W^{1,p}(\Omega)$ and
\begin{equation*}
\|w\|_{W^{1,p}(\Omega)}\leq C(\|w\|_{H^1(\Omega)}+\|\nabla u\|_{L^p(\Omega)})\leq C(\|u\|_{H^2(\Omega)}+\|u\|_{W^{1,p}(\Omega)})\leq C\|f\|_{H^{\frac{3}{2}}(\partial\Omega)}
\end{equation*}
where we use \eqref{est_u_h2}\eqref{est_u_w1p} in the last step, and $C$ also depends on $N$ defined in the admissible set $\mathcal{A}$.
Let $i$ vary from $1$ to $d$, we know $u\in W^{2,p}(\Omega)$ and
\begin{equation}
\|u\|_{W^{2,p}(\Omega)}\leq C\|f\|_{H^{\frac{3}{2}}(\partial\Omega)}.\label{est_u_w2p}
\end{equation}
\begin{lma}
Suppose $\sigma\in \mathcal{A}$, $\sigma+\delta\sigma\in \mathcal{A}$ with $\delta\sigma=0$ on $\partial\Omega$, and $f\in H^{\frac{3}{2}}(\partial\Omega)$. Let $u$ be the solution of \eqref{eq_int_u} and $v$ be the solution of
\begin{equation*}
\left\{\begin{array}{ll}\Div ((\sigma+\delta\sigma)\nabla v)=0&\quad \mbox{in}\
\Omega\\
v=f&\quad \mbox{on}\
\partial\Omega.\end{array} \right.
\end{equation*}
We have the following estimates
\begin{equation}
\|v-u\|_{H^1(\Omega)}
\leq C\| \delta\sigma\|^{1-\frac{2}{p}}_{L^2(\Omega)}\|f\|_{H^{\frac{1}{2}}(\partial \Omega)}\label{h1_cont}
\end{equation}

\begin{equation}
\|v-u\|_{H^2(\Omega)}
\leq C\| \delta\sigma\|^{1-\frac{2}{p}}_{H_0^1(\Omega)}\|f\|_{H^{\frac{3}{2}}(\partial \Omega)}\label{h2_cont}
\end{equation}
where $p>2$ from Theorem~\ref{thm_Meyers} and $C$ may depend on $c_1$, $c_2$, $N$, $\Omega$ and $p$.
\end{lma}
\noindent {\it{Proof:}} Clearly, $v-u$ satisfies
\begin{equation*}
\left\{\begin{array}{ll}\Div (\sigma\nabla (v-u))=-\Div (\delta\sigma\nabla v)&\quad \mbox{in}\
\Omega\\
v-u=0&\quad \mbox{on}\
\partial\Omega.\end{array} \right.
\end{equation*}
From the standard elliptic theory, we have
\begin{equation*}
\|v-u\|_{H^1(\Omega)}\leq C\|\delta\sigma\nabla v\|_{L^2(\Omega)}.
\end{equation*}
Applying the estimate \eqref{est_u_w1p} to $v$ and using the H\"{o}lder inequality, we obtain
\begin{equation*}
\| \delta\sigma\nabla v\|_{L^2(\Omega)}\leq \| \delta\sigma\|_{L^q(\Omega)}\| \nabla v\|_{L^p(\Omega)}\leq \| \delta\sigma\|_{L^q(\Omega)}\|f\|_{H^{\frac{1}{2}}(\partial \Omega)}
\end{equation*}
where $p>2$ from Theorem~\ref{thm_Meyers} and $q>2$ is such that $\frac{1}{p}+\frac{1}{q}=\frac{1}{2}$. From
\begin{equation*}
\| \delta\sigma\|_{L^q(\Omega)}=\left(\int_\Omega|\delta\sigma|^{q-2}|\delta\sigma|^2dx\right)^{\frac{1}{q}}\leq C\left(\int_\Omega|\delta\sigma|^2dx\right)^{\frac{1}{q}}=C\| \delta\sigma\|^{1-\frac{2}{p}}_{L^2(\Omega)},
\end{equation*}
we then get \eqref{h1_cont}.

From the standard elliptic theory, we also have
\begin{equation*}
\|v-u\|_{H^2(\Omega)}\leq C\|\Div (\delta\sigma\nabla v)\|_{L^2(\Omega)}
\leq C\|\nabla \delta\sigma\cdot\nabla v\|_{L^2(\Omega)}+C\| \delta\sigma\Delta v\|_{L^2(\Omega)}.
\end{equation*}
Similarly, applying the estimate \eqref{est_u_w1p} to $v$, we obtain
\begin{equation*}
\|\nabla \delta\sigma\cdot\nabla v\|_{L^2(\Omega)}\leq \| \nabla\delta\sigma\|_{L^q(\Omega)}\| \nabla v\|_{L^p(\Omega)}
\leq C\| \nabla\delta\sigma\|_{L^q(\Omega)}\|f\|_{H^{\frac{3}{2}}}\leq C\| \nabla\delta\sigma\|^{1-\frac{2}{p}}_{L^2(\Omega)}\|f\|_{H^{\frac{3}{2}}(\partial \Omega)},
\end{equation*}
and applying the estimate \eqref{est_u_w2p} to $v$, we obtain
\begin{equation*}
\| \delta\sigma\Delta v\|_{L^2(\Omega)}\leq \| \delta\sigma\|_{L^q(\Omega)}\| \Delta v\|_{L^p(\Omega)}
\leq C\| \delta\sigma\|_{L^q(\Omega)}\|f\|_{H^{\frac{3}{2}}(\partial \Omega)}\leq C\| \delta\sigma\|^{1-\frac{2}{p}}_{L^2(\Omega)}\|f\|_{H^{\frac{3}{2}}(\partial \Omega)}.
\end{equation*}
Hence  \eqref{h2_cont} holds.\hfill $\qed$


\begin{thm}
There exists at least one minimizer to the functional $R(\sigma)$.
\label{min_existence}
\end{thm}
\noindent {\it{Proof:}}  Since $R(\sigma)$ is nonnegative, there exists a minimizing sequence
$\{\sigma_n\}\subset\mathcal{A}$ such that 
\begin{equation*}
R(\sigma_n)\to R:=\liminf_{\sigma\in\mathcal{A}}R(\sigma)\quad \mbox{as}\quad  n\to\infty.
\end{equation*}
Clearly, $\sigma_n-\sigma^0$ is uniformly bounded in $H_0^1(\Omega)$. Thus, there exists a weakly convergent subsequence of $\{\sigma_n\}$, still denoted by $\{\sigma_n\}$, such that 
\begin{equation}
\sigma_n-\sigma^0\rightharpoonup\hat\sigma-\sigma^0 \ \mbox{ weakly in }\  H_0^1(\Omega) \ \mbox{ and } \|\hat\sigma-\sigma^0\|_{H_0^1(\Omega)}\leq \liminf_{n\to\infty} \|\sigma_n-\sigma^0\|_{H_0^1(\Omega)}.\label{weak_conv_1}
\end{equation}
From the compact Sobolev embedding $H_0^1(\Omega)\to L^2(\Omega)$, we have $\sigma_n\to\hat\sigma$ in $L^2(\Omega)$. Let $\hat u$ and $u_n$ be the solutions of \eqref{eq_int_u} with $\sigma=\hat\sigma$ and $\sigma=\sigma_n$ $(n=1, 2, \cdots)$. Applying \eqref{h1_cont} for 
$\sigma=\hat\sigma$ and $\delta\sigma=\sigma_n-\hat\sigma$, we have $\|u_n-\hat u\|_{H^1(\Omega)}\to 0$, that is, 
\begin{equation}
u_n\to \hat u \ \mbox{ in }\  H^1(\Omega).\label{hat_u_h1}
\end{equation}
Applying \eqref{h2_cont}, we then have 
\begin{eqnarray*}
\|u_n-\hat u\|_{H^2(\Omega)}
\leq C\| \sigma_n-\hat\sigma\|^{1-\frac{2}{p}}_{H_0^1(\Omega)}\|f\|_{H^{\frac{3}{2}}(\partial \Omega)}.
\end{eqnarray*}
Since $\sigma_n$ is uniformly bounded in $H^1(\Omega)$, from the above inequality, we know $u_n$ is uniformly bounded in $H^2(\Omega)$. Thus, there exists a weakly convergent subsequence of $\{u_n\}$, still denoted by $\{u_n\}$, such that 
\begin{equation*}
u_n\rightharpoonup\hat{\hat{u}} \ \mbox{ weakly in }\  H^2(\Omega).
\end{equation*}
From the compact Sobolev embedding $H^2(\Omega)\to H^1(\Omega)$, we have $u_n\to\hat{\hat{u}}$ in $H^1(\Omega)$. In view of \eqref{hat_u_h1}, we know $\hat u=\hat{\hat{u}}$. So 
\begin{equation*}
u_n\rightharpoonup\hat u\ \mbox{ weakly in }\  H^2(\Omega).
\end{equation*}
The trace operator mapping from $u$ to $\frac{\p u}{\p \nu}|_{\partial \Omega}$ is bounded from $H^2(\Omega)$ to $H^{\frac{1}{2}}(\partial \Omega)$. The embedding from $H^{\frac{1}{2}}(\partial \Omega)$ to $L^2(\partial \Omega)$ is compact, so we have
\begin{equation*}
\frac{\p u_n}{\p \nu}\to\frac{\p \hat u}{\p \nu}\ \mbox{ in }\  L^2(\partial \Omega)
\end{equation*}
and hence
\begin{equation}
\|F(\sigma_n,f)-F(\hat\sigma,f)\|_{L^2(\partial \Omega)}=\|\sigma_n\frac{\p u_n}{\p \nu}-\hat\sigma\frac{\p \hat u}{\p \nu}\|_{L^2(\partial \Omega)}=\|\hat\sigma (\frac{\p u_n}{\p \nu}-\frac{\p \hat u}{\p \nu})\|_{L^2(\partial \Omega)}\to 0.\label{F_cont_1}
\end{equation}
The existence of the the minimizer then follows from the continuity of $F$ and the weak lower semicontinuity of the norm. In fact, from \eqref{weak_conv_1}\eqref{F_cont_1}, we have $R(\hat\sigma)=R$.

\subsection{Gauss-Newton  method.} 

%

\xs{We first introduce some notations. Let $(\ ,\ )_{H^1(\Omega)}$,
$(\ ,\ )_\Omega=(\ ,\ )_{L^2(\Omega)}$, and $<\ ,\ >_{\partial \Omega}=<\ ,\ >_{L^2(\partial \Omega)}$ denote the inner
products on $H^1(\Omega)$, $L^2(\Omega)$, and $L^2(\partial \Omega)$. Let $Id$ be the identity operator on $L^2(\Omega)$.}

To find the minimizer $\sigma$ to $\displaystyle \min_{\sigma\in\mathcal{A}} \ R(\sigma)$, iterative methods are commonly used. In this work, 
we use the  well-known Gauss-Newton method. 
The iterative procedure reads
\begin{equation}
\sigma^{k+1}=\sigma^{k}+\delta\sigma,\quad k=0, 1, 2, \dots \label{eq_int_sigma_k}
\end{equation}
with $\delta\sigma$ solving
\begin{equation}
D^2R(\sigma^{k})\delta\sigma=-DR(\sigma^{k})\label{eq_int_delta_sigma_k}
\end{equation}
where $DR$ and $D^2R$ are the first derivative operator and second derivative operator of $R$, respectively.

Next we derive the formulas for $DR$ and $D^2R$ in order to find the update $\delta\sigma$ in \eqref{eq_int_delta_sigma_k}. 
We start with the derivative formulas related to the Dirichlet-to-Neumann operator $F(\sigma, \cdot)$. 
 The derivative formulas for Dirichlet-to-Neumann operator are similar to the derivative formulas for Neumann-to-Dirichlet operator. 
Let $DF$ be the derivative of $F$ with respect to $\sigma$, $(DF)^*$ and $(D^2F)^*$ be the adjoints of the first and second derivatives of $F$ with respect to $\sigma$.

\begin{lma}\label{lma:df}
Suppose $\sigma\in \mathcal{A}$, $\sigma+\delta\sigma\in \mathcal{A}$ with $\delta\sigma=0$ on $\partial\Omega$, and $f\in H^{\frac{3}{2}}(\partial\Omega)$. We have $DF(\sigma, f): H_0^1(\Omega)\to H^{\frac{1}{2}}(\partial\Omega)\subset L^2(\partial\Omega)$ is given by
\begin{equation}
DF(\sigma, f)(\delta\sigma)=\sigma\frac{\p \delta u}{\p \nu}\Big|_{\partial\Omega}\label{eq_int_DF}
\end{equation}
where $\delta u$ is the solution of 
\begin{equation}
\left\{\begin{array}{ll}\Div (\sigma\nabla \delta u)=-\Div(\delta\sigma\nabla u)&\quad \mbox{in}\
\Omega\\
\delta u=0&\quad \mbox{on}\
\partial\Omega\end{array} \right.
\label{eq_int_delta_u}
\end{equation}
with $u$ being the solution of \eqref{eq_int_u}. 
\end{lma}

\noindent{\bf Proof :} Let $\tilde{u}$ be the solution of
\begin{equation}
\left\{\begin{array}{ll}\Div ((\sigma+\varepsilon\delta\sigma)\nabla \tilde{u})=0&\quad \mbox{in}\
\Omega\\
\tilde{u}=f&\quad \mbox{on}\
\partial\Omega.\end{array} \right.
\label{eq_int_u_tilde}
\end{equation}
Then
\begin{equation*}
F(\sigma+\varepsilon\delta\sigma, f)=(\sigma+\varepsilon\delta\sigma)\frac{\p \tilde{u}}{\p \nu}\Big|_{\partial\Omega}= \sigma\frac{\p \tilde{u}}{\p \nu}\Big|_{\partial\Omega},
\end{equation*}
where we use $\delta\sigma=0$ on $\partial\Omega$ in the last step. \\
Denote 
\begin{equation*}
\delta u=\lim_{\varepsilon\to 0}\frac{\tilde{u}-u}{\varepsilon}.
\end{equation*}
Direct computation shows
\begin{eqnarray*}
&&DF(\sigma, f)(\delta\sigma)=\frac{d}{d\varepsilon}F(\sigma+\varepsilon\delta\sigma,f)\big |_{\varepsilon=0}=\lim_{\varepsilon\to 0}\frac{F(\sigma+\varepsilon\delta\sigma,f)-F(\sigma,f)}{\varepsilon}\\
&=&\lim_{\varepsilon\to 0}\frac{\sigma\frac{\p \tilde{u}}{\p \nu}-\sigma\frac{\p u}{\p \nu}}{\varepsilon}\Big|_{\partial\Omega}=\lim_{\varepsilon\to 0}\sigma\frac{\p}{\p \nu}(\frac{\tilde{u}-u}{\varepsilon})\Big|_{\partial\Omega}=\sigma\frac{\p \delta u}{\p \nu}\Big|_{\partial\Omega}.
\end{eqnarray*}
Next we show that $\delta u$ satisfies \eqref{eq_int_delta_u}. From \eqref{eq_int_u}\eqref{eq_int_u_tilde}, clearly $\delta u=0$ on $\partial \Omega$, and by taking the difference of these two equations in $\Omega$ we have
\begin{equation*}
\Div (\sigma\nabla (\tilde{u}-u))=-\varepsilon \Div(\delta\sigma\nabla \tilde{u}).
\end{equation*}
So
\begin{eqnarray*}
\Div (\sigma\nabla (\frac{\tilde{u}-u}{\varepsilon}))&=&- \Div(\delta\sigma\nabla \tilde{u})=-\Div(\delta\sigma\nabla (\tilde{u}-u))-\Div(\delta\sigma\nabla u)\\
&=&-\varepsilon\Div(\delta\sigma\nabla (\frac{\tilde{u}-u}{\varepsilon}))-\Div(\delta\sigma\nabla u).
\end{eqnarray*}
Letting $\varepsilon\to 0$ we know $\delta u$ satisfies \eqref{eq_int_delta_u}.

\begin{lma} \label{lma:df*}
Suppose $\sigma\in \mathcal{A}$ and $f\in H^{\frac{3}{2}}(\partial\Omega)$. We have $(DF)^*(\sigma, f): H^{\frac{1}{2}}(\partial\Omega)\subset L^2(\partial\Omega)\to H_0^1(\Omega)$ is given by
\begin{equation}
(DF)^*(\sigma, f)(\varphi)=w\label{eq_int_DF_star}
\end{equation} 
where $w$ is the solution of 
\begin{equation}\label{eq:sobolev}
\left\{\begin{array}{ll}-\Delta w+w=\nabla u \cdot \nabla u^*&\quad \mbox{in}\
\Omega\\
w=0&\quad \mbox{on}\
\partial\Omega\end{array} \right.
\end{equation}
with $u^*$ being the solution of 
\begin{equation}
\left\{\begin{array}{ll}\Div (\sigma\nabla u^*)=0&\quad \mbox{in}\
\Omega\\
u^*=\varphi &\quad \mbox{on}\
\partial\Omega\end{array} \right.
\label{eq_int_u_star}
\end{equation}
and $u$ being the solution of \eqref{eq_int_u}. 
\end{lma}

\noindent{\bf Proof:} For any $\delta\sigma\in H_0^1(\Omega)$, from \eqref{eq_int_DF} and the boundary condition in \eqref{eq_int_u_star}, we have 
\begin{equation}
(\delta\sigma, (DF)^*(\sigma,f)\varphi)_{H^1(\Omega)}=<DF(\sigma,f)\delta\sigma, \varphi>_{L^2(\partial\Omega)}=
\int_{\partial\Omega}\sigma\frac{\p \delta u}{\p \nu}u^*ds.\label{eq_int_df_star_p0}
\end{equation}
Next we show that
\begin{equation}
\int_{\partial\Omega}u^*\sigma\frac{\p \delta u}{\p \nu}ds=\int_{\Omega}\delta\sigma\nabla u\cdot\nabla u^*dx.\label{eq_int_df_star_p1}
\end{equation}
Multiplying $u^*$ to both sides of the equation \eqref{eq_int_delta_u} in $\Omega$ and integrating by parts, we have
\begin{equation}
\int_{\partial\Omega}u^*\sigma\frac{\p \delta u}{\p \nu}ds-\int_{\Omega}\sigma\nabla\delta u\cdot\nabla u^*dx=
-\int_{\partial\Omega}u^*\delta\sigma\frac{\p u}{\p \nu}ds+\int_{\Omega}\delta\sigma\nabla u\cdot\nabla u^*dx.\label{eq_int_df_star_p1_2}
\end{equation}
Since $\delta \sigma=0$ on $\partial\Omega$, the first term on the right hand side of \eqref{eq_int_df_star_p1_2} is $0$. Multiplying $\delta u$ to both sides of the equation \eqref{eq_int_u_star} in $\Omega$ and integrating by parts, we have
\begin{equation*}
\int_{\partial\Omega}\delta u\sigma\frac{\p u^*}{\p \nu}ds-\int_{\Omega}\sigma\nabla u^*\cdot\nabla \delta udx=0.
\end{equation*}
In view of $\delta u=0$ in $\partial\Omega$, we know that $\int_{\Omega}\sigma\nabla\delta u\cdot\nabla u^*dx=0$, that is, the second term on the left hand side of \eqref{eq_int_df_star_p1_2} is also $0$. So \eqref{eq_int_df_star_p1} holds. 

From \eqref{eq_int_df_star_p0}\eqref{eq_int_df_star_p1}, we have
\begin{equation*}
(\delta\sigma, w)_{H^1(\Omega)}=(\delta\sigma, (DF)^*(\sigma,f)\varphi)_{H^1(\Omega)}=(\delta\sigma, \nabla u \cdot \nabla u^*)_{L^2(\Omega)},
\end{equation*}
that is,
\begin{equation*}
\int_{\Omega}(\delta\sigma w+\nabla \delta\sigma\cdot \nabla w) dx=\int_{\Omega}\delta\sigma \nabla u \cdot \nabla u^*dx.
\end{equation*}
Hence $w$ is the solution of \eqref{eq:sobolev}.



\begin{rem} $w$ is known as the Sobolev gradient (see, for example \cite{Ne, JM}). The regularity requirement for the conductivity in the admissible set is also used here. 
\end{rem}

After getting the derivatives formulas related to $F$, we now study the formulas for $DR$ and $D^2R$. 
\begin{lma} \label{lemma}
The first and second derivatives of $R$ of \eqref{eq_int_R} are 
\begin{equation}
DR(\sigma)=\sum_{j=1}^{M}(DF)^*(\sigma,f_j)(F(\sigma,f_j)-g_j^\delta)+\alpha(Id-\Delta)(\sigma-\sigma^0)
\label{eq_int_DR}
\end{equation}
and
\begin{equation}
D^2R(\sigma)=\sum_{j=1}^{M}\left[(DF)^*(\sigma,f_j)DF(\sigma,f_j)+(D^2F)^*(\sigma,f_j)(F(\sigma,f_j)-g_j^\delta)\right]+\alpha (Id-\Delta)\,,
\label{eq_int_D2R}
\end{equation}
respectively.
\end{lma}
\noindent{\bf Proof :} We directly compute the derivative of $R$ at $\sigma$ in the direction $\delta\sigma$. From \eqref{eq_int_R} we know
\begin{eqnarray*}
&R(\sigma+\varepsilon\delta\sigma)=&\frac{1}{2}\sum_{j=1}^{M}<F(\sigma+\varepsilon\delta\sigma,f_j)-g_j^\delta, F(\sigma+\varepsilon\delta\sigma,f_j)-g_j^\delta>_{\partial \Omega}\\
&&+\frac{\alpha}{2}(\sigma+\varepsilon\delta\sigma-\sigma^0, \sigma+\varepsilon\delta\sigma-\sigma^0)_{\Omega}
+\frac{\alpha}{2}(\nabla(\sigma+\varepsilon\delta\sigma-\sigma^0), \nabla(\sigma+\varepsilon\delta\sigma-\sigma^0))_{\Omega}.
\end{eqnarray*}
So
\begin{eqnarray}
&&DR(\sigma)(\delta\sigma)=\frac{d}{d\varepsilon}R(\sigma+\varepsilon\delta\sigma)\big |_{\varepsilon=0}\nonumber\\
&=&\frac{1}{2}\sum_{j=1}^{M}\Big(<DF(\sigma,f_j)\delta\sigma,F(\sigma,f_j)-g_j^\delta>_{\partial \Omega}+<F(\sigma,f_j)-g_j^\delta,DF(\sigma,f_j)\delta\sigma>_{\partial \Omega}\Big)\nonumber\\
&& +\frac{\alpha}{2}\Big((\delta\sigma,\sigma-\sigma^0)_{\Omega}+(\sigma-\sigma^0,\delta\sigma)_{\Omega}\Big)+\frac{\alpha}{2}\Big((\nabla\delta\sigma,\nabla(\sigma-\sigma^0))_{\Omega}+(\nabla(\sigma-\sigma^0),\nabla\delta\sigma)_{\Omega}\Big)\nonumber\\
&=&\sum_{j=1}^{M}<F(\sigma,f_j)-g_j^\delta,DF(\sigma,f_j)\delta\sigma>_{\partial \Omega}+\alpha(\sigma-\sigma^0,\delta\sigma)_{\Omega}+\alpha(\nabla(\sigma-\sigma^0),\nabla\delta\sigma)_{\Omega}\label{eq_int_DR_last2}\\
&=&\sum_{j=1}^{M}((DF)^*(\sigma,f_j)(F(\sigma,f_j)-g_j^\delta),\delta\sigma)_{\Omega}+\alpha(\sigma-\sigma^0,\delta\sigma)_{\Omega}+\alpha(-\Delta(\sigma-\sigma^0),\delta\sigma)_{\Omega},\nonumber
\end{eqnarray}
where we use the integration by part for the last term in \eqref{eq_int_DR_last2} and $\delta\sigma=0$ on the boundary. So \eqref{eq_int_DR} holds.

We then compute the bilinear second derivative of $R$ at $\sigma$ in the direction $\delta\sigma$. From \eqref{eq_int_DR_last2} we know
\begin{eqnarray*}
DR(\sigma+\varepsilon\delta\sigma)(\delta\sigma)&=&\sum_{j=1}^{M}<F(\sigma+\varepsilon\delta\sigma,f_j)-g_j^\delta,DF(\sigma+\varepsilon\delta\sigma,f_j)\delta\sigma>_{\partial \Omega}\\
&&+\alpha(\sigma+\varepsilon\delta\sigma-\sigma^0,\delta\sigma)_{\Omega}+\alpha(\nabla(\sigma+\varepsilon\delta\sigma-\sigma^0),\nabla\delta\sigma)_{\Omega}.
\end{eqnarray*}
So
\begin{eqnarray*}
&&D^2R(\sigma)(\delta\sigma, \delta\sigma)=\frac{d}{d\varepsilon}DR(\sigma+\varepsilon\delta\sigma)\delta\sigma |_{\varepsilon=0}\\
&=&\sum_{j=1}^{M}\Big(<DF(\sigma,f_j)\delta\sigma, DF(\sigma,f_j)\delta\sigma>_{\partial \Omega}+<F(\sigma,f_j)-g_j^\delta,D^2F(\sigma,f_j)(\delta\sigma, \delta\sigma)>_{\partial \Omega}\Big)\\
&&\quad +\alpha(\delta\sigma,\delta\sigma)_{\Omega}+\alpha(\nabla\delta\sigma,\nabla\delta\sigma)_{\Omega}\\
&=&\sum_{j=1}^{M}\Big(<(DF)^*(\sigma,f_j)DF(\sigma,f_j)\delta\sigma, \delta\sigma>_{\partial \Omega}+<(D^2F)^*(\sigma,f_j)(F(\sigma,f_j)-g_j^\delta, \delta\sigma), \delta\sigma>_{\partial \Omega}\Big)\\
&&\quad +\alpha(\delta\sigma,\delta\sigma)_{\Omega}+\alpha(-\Delta\delta\sigma,\delta\sigma)_{\Omega},
\end{eqnarray*}
which proves \eqref{eq_int_D2R}. \hfill $\qed$

\vspace{.5\baselineskip}
Now we can apply the formulas \eqref{eq_int_DR}, \eqref{eq_int_D2R} to compute $\delta\sigma$ in \eqref{eq_int_delta_sigma_k}. For simplicity, we 
also ignore the term $(D^2F)^*$ in $D^2R$ and solve the following linear equation without second derivative
\begin{eqnarray}
&&\left(\sum_{j=1}^{M}(DF)^*(\sigma^{k},f_j)DF(\sigma^{k},f_j)+\alpha(Id-\Delta)\right)\delta\sigma\nonumber\\
&=&-\sum_{j=1}^{M}(DF)^*(\sigma^{k},f_j)(F(\sigma^{k},f_j)-g_j^\delta)-\alpha(Id-\Delta)(\sigma^{k}-\sigma^0).\label{eq_int_omega}
\end{eqnarray}
We will use the conjugate gradient method (see, for example \cite{EHN}) to solve \eqref{eq_int_omega}. Typically this method only needs a small number of iteration steps by generating orthogonal residuals. The conjugate gradient methods for EIT related problems are studied in \cite{Ha, LR1, LR2} and the references therein.

\section{The MD-LDG method}

\xs{In this section, we introduce the numerical methods to solve the derivative operators $DF$ and $(DF)^*$ in
\eqref{eq_int_DF} and \eqref{eq_int_DF_star}, which involve several elliptic type equations.}
The minimal dissipation local discontinuous Galerkin method (MD-LDG) \cite{CoDo07} is used to solve
all the partial differential equations in each iteration step.
MD-LDG method is a special LDG
method for which the stabilization parameters are taken to be identically
zero on all interior faces.

Since our numerical examples are in two dimensions, let us illustrate the MD-LDG formulation on the model problem on the domain $\Omega\in\mathbb{R}^2$.
We remark that the formulation of MD-LDG in higher dimensions is similar.
\begin{equation}\label{eq:modeldg}
\left\{\begin{array}{ll}
\Div  (\sigma(x,y)\nabla u)=-r(x,y)&\quad \mbox{in}\
\Omega\\
u=b(x,y)&\quad \mbox{on}\
\partial\Omega\end{array} \right.
\end{equation}
where $\sigma(x,y)\in L^{\infty}(\Omega)$ satisfying $0<c_1<\sigma(x,y)<c_2$, $r(x,y)\in L^2(\Omega)$ and $b(x,y)\in L^2(\partial\Omega)$.

In order to define the LDG method, we rewrite \eqref{eq:modeldg} into a system of the first order
equations
\begin{equation}\label{eq:dg}
\begin{array}{ll}
\sigma(x,y)^{-1}{\bf q}=\nabla u\\
\nabla \cdot {\bf q}=-r(x,y)
\end{array} 
\end{equation}
where ${\bf q}=(q_1,q_2)$ is a vector function. Then we introduce the finite element spaces associated to the triangulation
$\Omega_h=\{K\}$ of $\Omega$ of shape-regular tetrahedra $K$.
We set
\[V_h=\{v\in L^2(\Omega): v|_K\in P^k(K), \forall K\in \Omega_h\}\]
\[{\bf W}_h=\{{\bf w}\in L^2(\Omega): {\bf w}|_K\in [P^k(K)]^2, \forall K\in \Omega_h\}\]
where $P^k(K)$ denotes the set of all polynomials of degree at most $k$ on $K$.
LDG method is to find $u_h \in V_h$ and ${\bf q}_h \in {\bf W}_h$ such that
for all $K\in\Omega_h$ and all test functions $v \in V_h$ and ${\bf w} \in {\bf W}_h$ we have 
\begin{equation}
\int_K \sigma(x,y)^{-1}{\bf q}_h\cdot {\bf w} \; dxdy+\int_K u_h \nabla \cdot {\bf w}\; dxdy-\int_{\partial K}{\widehat u_h} {\bf n}_K\cdot{\bf w}\;dS=0
\end{equation}

\begin{equation}
\int_K {\bf q}_h\cdot \nabla v  \; dxdy-\int_{\partial K}{\widehat{\bf q}_h}\cdot {\bf n}_K v\;dS-\int_K r(x,y)v\;dxdy=0
\end{equation}
where ${\bf n}_K$ is the outward normal unit vector to the $\partial K$. 

Next we define the  numerical fluxes ${\widehat u_h}$ and ${\widehat{\bf q}_h}$.
For a scalar valued function $u$, we define the average $\lb u \rb$ and the 
jump $\lbk u \rbk$ as follows. Let $e$ be an interior edge shared by 
elements $K_1$ and $K_2$. Define the unit normal vectors $\bld{n}_1$ and 
$\bld{n}_2$ on $e$ pointing exterior to $K_1$ and $K_2$, respectively. 
With $u_i:=u|_{\partial K_i}$, we set
\begin{equation*}
\lb u \rb=\frac{1}{2}(u_1+u_2),\quad \lbk u \rbk=u_1 \bld{n}_1+u_2 
\bld{n}_2.
\end{equation*}
For a vector-valued 
function $\bq$, we define $\bq_1$ and $\bq_2$ analogously and set
\begin{equation*}
\lb \bq \rb=\frac{1}{2}(\bq_1+\bq_2),\quad \lbk \bq \rbk=\bq_1 \cdot 
\bld{n}_1+\bq_2 \cdot \bld{n}_2.
\end{equation*}
We do not require either of the 
quantities $\lb u \rb$ or $\lbk \bq \rbk$ on boundary edges, and we leave 
them undefined. The fluxes are chosen as follows:
\begin{equation*}
\begin{array}{ll}
{\widehat u_h}=\lb u_h\rb+ {\boldsymbol\beta} \cdot \lbk u_h \rbk, & \mbox{on}\
\Gamma^0, \\
{\widehat u_h}=b, & \mbox{on}\
\partial\Omega, \\
\end{array} 
\end{equation*}
and
\begin{equation*}
\begin{array}{ll}
\widehat{\bf q}_h= \lb{\bf q}_h\rb-{\boldsymbol \beta} \lbk{\bf q}_h\rbk, & \mbox{on}\
\Gamma^0, \\
\widehat{\bf q}_h={\bf q}_h, & \mbox{on}\
\partial\Omega\cap\Gamma^- , \\
\widehat{\bf q}_h={\bf q}_h-\alpha(u_h-b){\bf n}, & \mbox{on}\
\partial\Omega\cap\Gamma^+, \\
\end{array} 
\end{equation*}
where ${\boldsymbol \beta} \cdot {\bf n}_K(e)=\frac{1}{2}\mbox{sign}({\bf v}\cdot {\bf n}_K(e))$ and
$\bf v$  is any nonzero piecewise constant vector.
 $\Gamma$ denotes the union of the boundaries of the element $K$ of $\Omega_h$ and
 $\Gamma^0$ denotes the interior boundaries $\Gamma^0:=\Gamma\backslash \partial \Omega$.
\[\Gamma^-=\{e\in \Gamma: {\bf v}\cdot {\bf n}_e<0\}, \quad 
\Gamma^+=\Gamma\backslash\Gamma^-.\]
The stabilization parameter $\alpha$ is chosen as $\mathcal O(1/h)$.

We refer the error estimate results and proofs in \cite{CoDo07}.
\begin{thm}\label{thm}
Suppose that  $\Omega$ is convex and that the exact solution $({\bf q}, u)$ of \eqref{eq:modeldg} belongs to ${\bf H}^r(\Omega_h)\times H^{r+1}(\Omega_h)$, for some $r\in[1,k]$.
Let $({\bf q}_h, u_h)\in {\bf W}_h\times V_h $ be the approximated solution by MD-LDG defined above, then we have
\begin{equation}
||{\bf q}-{\bf q}_h||_{L^2(\Omega_h)}\le C_1({\bf q}, u) h^r,
\end{equation}

\begin{equation}
||u-u_h||_{L^2(\Omega_h)}\le C_2({\bf q}, u) h^{r+1},
\end{equation}
where $C_1$ and $C_2$ are dependent of $\bf q$ and $u$ but independent of $h$.
\end{thm}
The order of convergence of the solution $u$ by LDG with $P^k$ polynomial space is order $k+1$ which is optimal. The order of convergence of  $\bf q$ is of order $k$ ( except in 1D, it is of order $k+1$).

\begin{rem}
There are different ways of defining {\bf q} in LDG method in \eqref{eq:dg}. Our definition is natural because of the Dirichlet-to-Neumann map
in \eqref{eq:f}.
\end{rem}

\section{Numerical  algorithms}
In this section we precisely describe our numerical algorithms. The Gauss-Newton method is used to find the minimizer of \eqref{eq_int_R}. 
The iteration reads \eqref{eq_int_sigma_k}, which is the outer iteration. In each iteration, from the analysis in Section~2, we need to solve \eqref{eq_int_omega}. It will be solved by the conjugate gradient algorithm, which is the inner iteration.

\subsection{The Gauss-Newton algorithm} We describe the initialization, stopping criterion and the iteration steps for \eqref{eq_int_sigma_k}.

\vspace{.5\baselineskip}

\noindent {\bf Initialization:}

Given an initial guess for conductivity $\sigma^0$. Given $M$ measurements of voltage on the boundary $f_j$, $j=1, \dots, M$.
The exact current flux $g_j^{true}=\sigma^{true}\frac{\partial u_j}{\partial \nu}|_{\partial \Omega}$ are precomputed   from 
\begin{equation}
\left\{\begin{array}{ll}
\Div  (\sigma^{true}\nabla u_j)=0&\quad \mbox{in}\
\Omega\\
u_j=f_j&\quad \mbox{on}\
\partial\Omega\end{array} \right.
\end{equation}
by MD-LDG on a fine mesh.

We add the noise to the exact current flux in the following way
\begin{equation}
g_j^\delta=g_j^{true}+\varepsilon |g_j^{true}|\xi_j, \quad j=1,\dots,M
\end{equation}
where  $\xi_j$ follow the standard normal distribution.

\vspace{.5\baselineskip}

\noindent {\bf Stopping criterion:}

The iteration will be stopped when the error between the computed data and the measured data reaches the noisy level. \xs{More precisely, 
let $\delta$ be the $L^2$ norm of the noise level on the boundary 
\begin{equation*}
\sum_{j=1}^{M}\|g_j^\delta-g_j^{true}\|_{L^2(\partial\Omega)}\le \delta.
\end{equation*}}
\xs{We take
\begin{equation*}
\delta=\sum_{j=1}^{M}\|g_j^\delta-g_j^{true}\|_{L^2(\partial\Omega)}=\sum_{j=1}^{M}\| \varepsilon g_j^{true}\xi_j\|_{L^2(\partial\Omega)}
=\varepsilon\sum_{j=1}^{M}\|  g_j^{true}\xi_j\|_{L^2(\partial\Omega)}.
\end{equation*}}The iteration will be stopped when $\sum_{j=1}^{M}\|F(\sigma^k,f_j)-g_j^\delta\|_{L^2(\partial \Omega)}$ reaches the order of $\delta$. We choose 
$\tau>1$ and stop the iteration at the first occurrence of $k$ such that 
\begin{equation*}
\sum_{j=1}^{M}\|F(\sigma^k,f_j)-g_j^\delta\|_{L^2(\partial \Omega)}\leq\tau\delta.
\end{equation*}
We will discuss $\tau$ in the next subsection.

\vspace{.5\baselineskip}

\noindent {\bf Algorithm:}



Set $k=0$. Input a constant $\tau>1$ and a maximum number of iterations MaxOut.  Start the  iteration to solve for $\sigma$.

While $\Big(\sum_{j=1}^{M}\|F(\sigma^k,f_j)-g_j^\delta\|_{L^2(\partial\Omega)}>\tau\delta \mbox{ and } k<\mbox{MaxOut}\Big)$ do
\begin{enumerate}
\item
Solve the following linear equation for $\delta\sigma$ using the conjugate gradient method:
\begin{eqnarray}\label{eq:CG}
&&\left(\sum_{j=1}^{M}(DF)^*(\sigma^{k},f_j)DF(\sigma^{k},f_j)+\alpha(Id-\Delta) \right)\delta\sigma\nonumber\\
&=&-\sum_{j=1}^{M}(DF)^*(\sigma^{k},f_j)(F(\sigma^{k},f_j)-g_j^\delta)-\alpha(Id-\Delta)(\sigma^{k}-\sigma^0).
\end{eqnarray}


\item Set $\sigma^{k+1}=\sigma^{k}+\delta\sigma$.

\item Set $k:=k+1$.

\end{enumerate}

We obtain $\sigma=\sigma^k.$ 

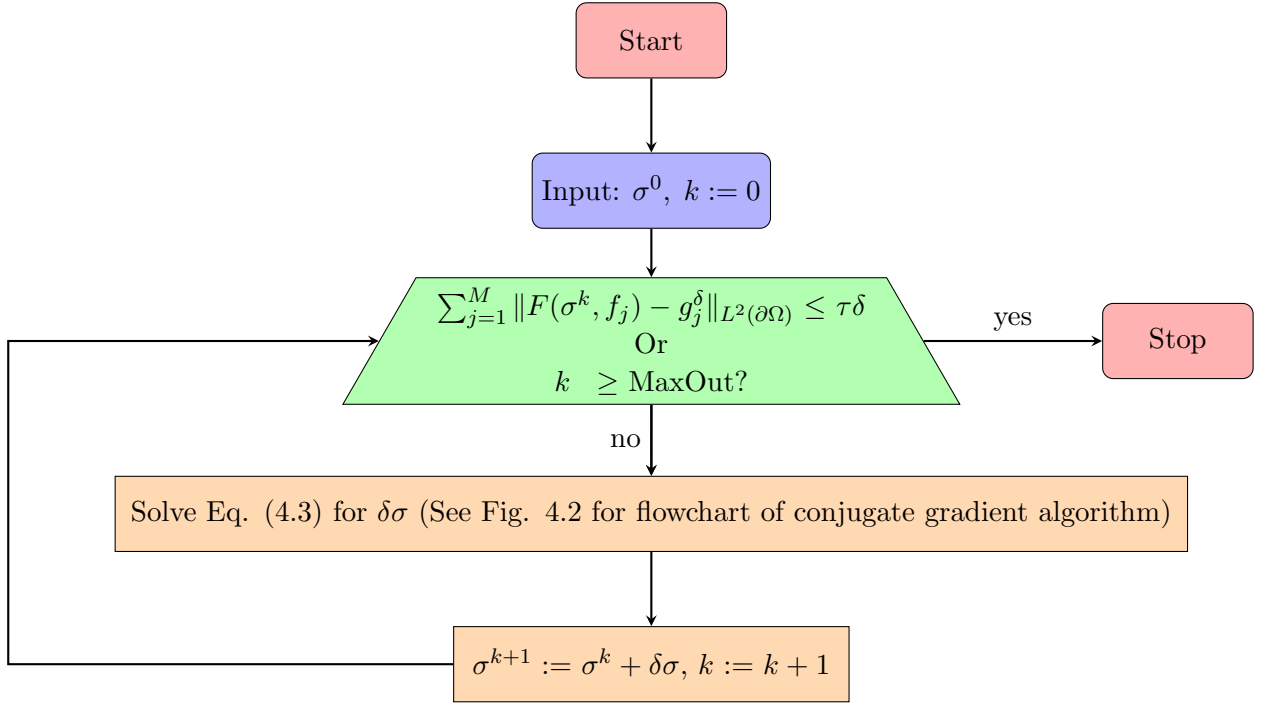
\begin{figure}

\begin{tikzpicture}[node distance=2cm,
     arr/.style = {draw=black!50, thick, -Latex}
]

\node (start) [startstop] {Start};
\node (in1) [io, below of=start] {Input: $\sigma^0,\; k:=0$};

\node (dec1) [decision, below of=in1] {$\sum_{j=1}^{M}\|F(\sigma^k,f_j)-g_j^\delta\|_{L^2(\partial\Omega)}\le\tau\delta$ \\
Or \\
$k\ge$ MaxOut?};
\node (pro1) [process2, below of=dec1, yshift=-0.3cm] {Solve Eq. (4.3) for $\delta\sigma$ (See Fig. \ref{flow:2} for flowchart of conjugate gradient algorithm)};
\node (pro2) [process, below of=pro1] {$\sigma^{k+1}:=\sigma^{k}+\delta\sigma$,
$k:=k+1$};

\node (stop) [startstop, right of=dec1, xshift=5cm] {Stop};

\draw [arrow] (start) -- (in1);
\draw [arrow] (in1) -- (dec1);
\draw [arrow] (dec1) -- (pro1);
\draw [arrow] (dec1) -- node[anchor=south] {yes} (stop);
\draw [arrow] (dec1) -- node[anchor=east] {no} (pro1);
\draw [arrow] (pro1) --(pro2);
\coordinate[above left=of pro1] (aux);
\draw [arrow] (pro2) -|(aux)|-(dec1);
\end{tikzpicture}

{\caption{Flowchart of Gauss-Newton algorithm}}
\label{flow:1}

\end{figure}

{A flowchart of Gauss-Newton  algorithm is shown in Figure \ref{flow:1}.}

\subsection{The conjugate gradient algorithm}
We use the conjugate gradient method to solve the linear equation \eqref{eq:CG}.

\vspace{.5\baselineskip}

\noindent {\bf Initialization:}

Denote the right hand side of \eqref{eq:CG} as
\[r_0=-\sum_{j=1}^{M}{(DF)^*}(\sigma^{k},f_j)(F(\sigma^{k},f_j)-g_j^\delta)-\alpha(Id-\Delta)(\sigma^{k}-\sigma^0),\]
where $(DF)^*$ is defined in Lemma \ref{lma:df*}, and \eqref{eq_int_DF_star} in Lemma \ref{lma:df*} is solved by
MD-LDG.

Given an initial guess $(\delta\sigma)_0=0$.
Set the initial direction \(p_0=r_0\).

\vspace{.5\baselineskip}

\noindent {\bf Stopping criterion:} The iteration will be stopped when the relative residual is smaller than a given tolerance $\rho$ $(0<\rho< 1)$. More precisely, we stop the iteration at the first occurrence of $l$ such that
\begin{equation*}
\sum_{j=1}^{M}||g_j^\delta-F(\sigma^{k},f_j)-DF(\sigma^{k},f_j)(\delta\sigma)_l||_{L^2(\partial\Omega)}<\rho \sum_{j=1}^{M} ||g_j^\delta-F(\sigma^{k},f_j)||_{L^2(\partial\Omega)}.
\end{equation*}
We also require $\rho^2\tau>2$ as in \cite{Ha}. In this paper, we fix $\tau=3 $ and $\rho= 0.9$. We would like to point out that it is not our purpose to choose the optimal numbers for $\rho$ and $\tau$. 

\vspace{.5\baselineskip}

\noindent {\bf Algorithm:}

Set $l=0$.  Input a constant $\rho$ and a maximum number of iterations MaxInn.  Start the  iteration to solve for $\delta\sigma$. 

While $\Big(\sum_{j=1}^{M}||g_j^\delta-F(\sigma^{k},f_j)-DF(\sigma^{k},f_j)(\delta\sigma)_l||_{L^2(\partial\Omega)}\geq\rho \sum_{j=1}^{M} ||g_j^\delta-F(\sigma^{k},f_j)||_{L^2(\partial\Omega)}$ and $l<MaxInn \Big)$ do

\begin{enumerate}

\item Set\[\alpha_l=\frac{||r_{l-1}||_{L^2(\Omega)}^2}{\sum_{j=1}^{M}||DF(\sigma^{k},f_j) p_{l-1}||_{L^2(\Omega)}^2+\alpha ||p_{l-1}||_{L^2(\Omega)}^2+\alpha ||\nabla p_{l-1}||_{L^2(\Omega)}^2},\]
where $DF$ is defined in Lemma \ref{lma:df}, and \eqref{eq_int_DF} in Lemma \ref{lma:df} is solved by
MD-LDG.

\item Set \((\delta\sigma)_l=(\delta\sigma)_{l-1}+\alpha_l p_{l-1}.\)

\item Set \[r_l=r_{l-1}-\alpha_l \left(\sum_{j=1}^{M}{(DF)^*}(\sigma^{k},f_j)DF(\sigma^{k},f_j)+\alpha(Id-\Delta) \right)p_{l-1}.\]

\item  Set \(\displaystyle\beta_l=\frac{||r_{l}||_{L^2(\Omega)}^2}{||r_{l-1}||_{L^2(\Omega)}^2}.\)

\item Set \(p_l=r_l+\beta_l p_{l-1}.\)

\item Set \(l:=l+1\).
\end{enumerate}

We obtain $\delta\sigma=(\delta\sigma)_l$. 

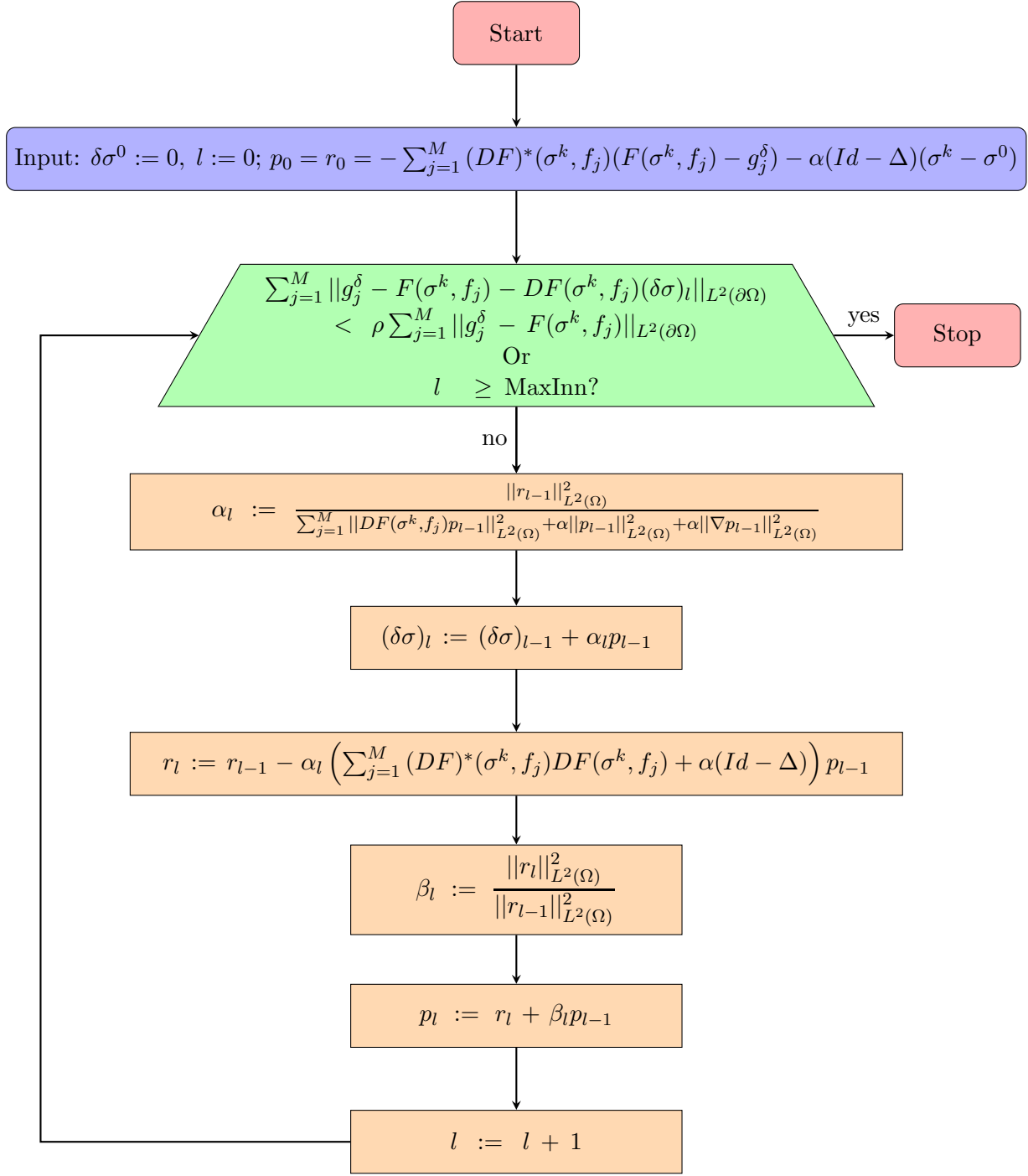
\begin{figure}

\begin{tikzpicture}[node distance=2cm]

\node (start) [startstop] {Start};
\node (in1) [io, below of=start] {Input: $\delta\sigma^0:=0,\; l:=0$; 
$p_0=r_0=-\sum_{j=1}^{M}{(DF)^*}(\sigma^{k},f_j)(F(\sigma^{k},f_j)-g_j^\delta)-\alpha(Id-\Delta)(\sigma^{k}-\sigma^0)$};

\node (dec1) [decision2, below of=in1,yshift=-0.8cm] {$\sum_{j=1}^{M}||g_j^\delta-F(\sigma^{k},f_j)-DF(\sigma^{k},f_j)(\delta\sigma)_l||_{L^2(\partial\Omega)}$

$<\rho \sum_{j=1}^{M} ||g_j^\delta-F(\sigma^{k},f_j)||_{L^2(\partial\Omega)}$

Or 

$l\ge\,$ MaxInn?};
\node (pro1) [process3, below of=dec1, yshift=-0.8cm] {$\alpha_l:=\frac{||r_{l-1}||_{L^2(\Omega)}^2}{\sum_{j=1}^{M}||DF(\sigma^{k},f_j) p_{l-1}||_{L^2(\Omega)}^2+\alpha ||p_{l-1}||_{L^2(\Omega)}^2+\alpha ||\nabla p_{l-1}||_{L^2(\Omega)}^2}$};
\node (pro2) [process, below of=pro1] {\((\delta\sigma)_l:=(\delta\sigma)_{l-1}+\alpha_l p_{l-1}\)};
\node (pro3) [process3, below of=pro2] {$r_l:=r_{l-1}-\alpha_l \left(\sum_{j=1}^{M}{(DF)^*}(\sigma^{k},f_j)DF(\sigma^{k},f_j)+\alpha(Id-\Delta) \right)p_{l-1}$};

\node (pro4) [process, below of=pro3] {$\displaystyle\beta_l:=\frac{||r_{l}||_{L^2(\Omega)}^2}{||r_{l-1}||_{L^2(\Omega)}^2}$};

\node (pro5) [process, below of=pro4] {\(p_l:=r_l+\beta_l p_{l-1}\)};
\node (pro6) [process, below of=pro5] {$l:=l+1$};

\node (stop) [startstop, right of=dec1, xshift=5cm] {Stop};

\draw [arrow] (start) -- (in1);
\draw [arrow] (in1) -- (dec1);
\draw [arrow] (dec1) -- (pro1);
\draw [arrow] (dec1) -- node[anchor=south] {yes} (stop);
\draw [arrow] (dec1) -- node[anchor=east] {no} (pro1);
\draw [arrow] (pro1) --(pro2);
\draw [arrow] (pro2) --(pro3);
\draw [arrow] (pro3) --(pro4);
\draw [arrow] (pro4) --(pro5);
\draw [arrow] (pro5) --(pro6);
\coordinate[above left=of pro3] (aux);
\draw [arrow] (pro6) -|(aux)|-(dec1);

\end{tikzpicture}

{\caption{Flowchart of conjugate gradient  algorithm}}
\label{flow:2}

\end{figure}

{A flowchart of conjugate gradient  algorithm is shown in Figure \ref{flow:2}.}

%

\section{Numerical results}
In this section, we will present several numerical experiments to demonstrate the 
 performance of the proposed numerical reconstruction method.
We first test our MD-LDG method for the forward problem.
Then we apply MD-LDG method as the forward solver to solve the iterative inverse problem. 
In  the numerical reconstructions, we use the following 4 measurements
\begin{equation}
f_1=\sin(x+y),\quad f_2=\cos(x+y),\quad f_3=\sin 2(x+y),\quad f_4=\cos 2(x+y).
\end{equation}
It is natural to choose the linearly independent sine and cosine functions as the measurements functions $f_j$. 
Note that more measurements may produce better results, but more computational cost.


\subsection{Example \ref{ex:1}: Convergence of forward problem}
\label{ex:1}

In the first example, we would like to test the convergence of our MD-LDG as the forward solver for \eqref{eq:modeldg}.
The convergence of MD-LDG for \eqref{eq:modeldg} is well known in the literature \xs{(see, for example, \cite{CoDo07})}.
We choose the exact solution 
$u=\sin(x+y)$ and the coefficient $\sigma=e^{-(x^2+y^2)}$.
The computational domain is a square $[0,1]\times[0,1]$.
The right hand side $r(x,y)$ and the boundary $b(x,y)$ in \eqref{eq:modeldg} are provided from the calculation of $u$.
We use the MD-LDG with $P^2$ polynomial space.
Table \ref{table:2} showed the $L^2$-errors and orders of accuracy of 
$u$, $\sigma  \frac{\partial u}{\partial \nu}$ (in this example $\sigma u_x=\sigma u_y$). We can see third order convergence for $u$ and second order for $\sigma  \frac{\partial u}{\partial \nu}$.
This is confirmed with the optimal convergence for $u$ and suboptimal convergence for $\bf q$  in Theorem \ref{thm}.

      \begin{table}
     \centering
   \caption{Example \ref{ex:1}: $L^2$-errors and orders of accuracy of MD-LDG $P^2$. }
    \smallskip
    \begin{tabular}{ccccc}\\\hline\hline
&$u$ && $\sigma  \frac{\partial u}{\partial \nu}$& \\\hline
$N$ & error & order  & error & order  \\\hline\hline
8$\times$ 8& 2.66E-05  &--  &1.23E-04    &--    \\

16$\times$16&3.22E-06  & 3.05   &1.82E-05   & 2.75    \\
32$\times$32&3.98E-07  &3.01    &  2.86E-06 &2.68      \\
64$\times$64&4.94E-08  &  3.01  &4.70E-07   &2.60      \\

\hline\hline
    \end{tabular}\label{table:1}
    \end{table}

\subsection{Example \ref{ex:2}: Reconstruction of EIT: one smooth blob}
\label{ex:2}

We consider a 2D problem  on the domain $\Omega=[-1,1]\times [-1,1]$. The true conductivity is given by
$\displaystyle \sigma(x)=\sigma_0(x)+e^{-8(x^2+(y-0.55)^2)}$ with the background
conductivity $\sigma_0=1$, same as  \cite{Zou17}. \xs{We take the background
conductivity as our initial guess.}

\begin{figure}
 \begin{center}
    \includegraphics[width = 2in,height=1.8in] {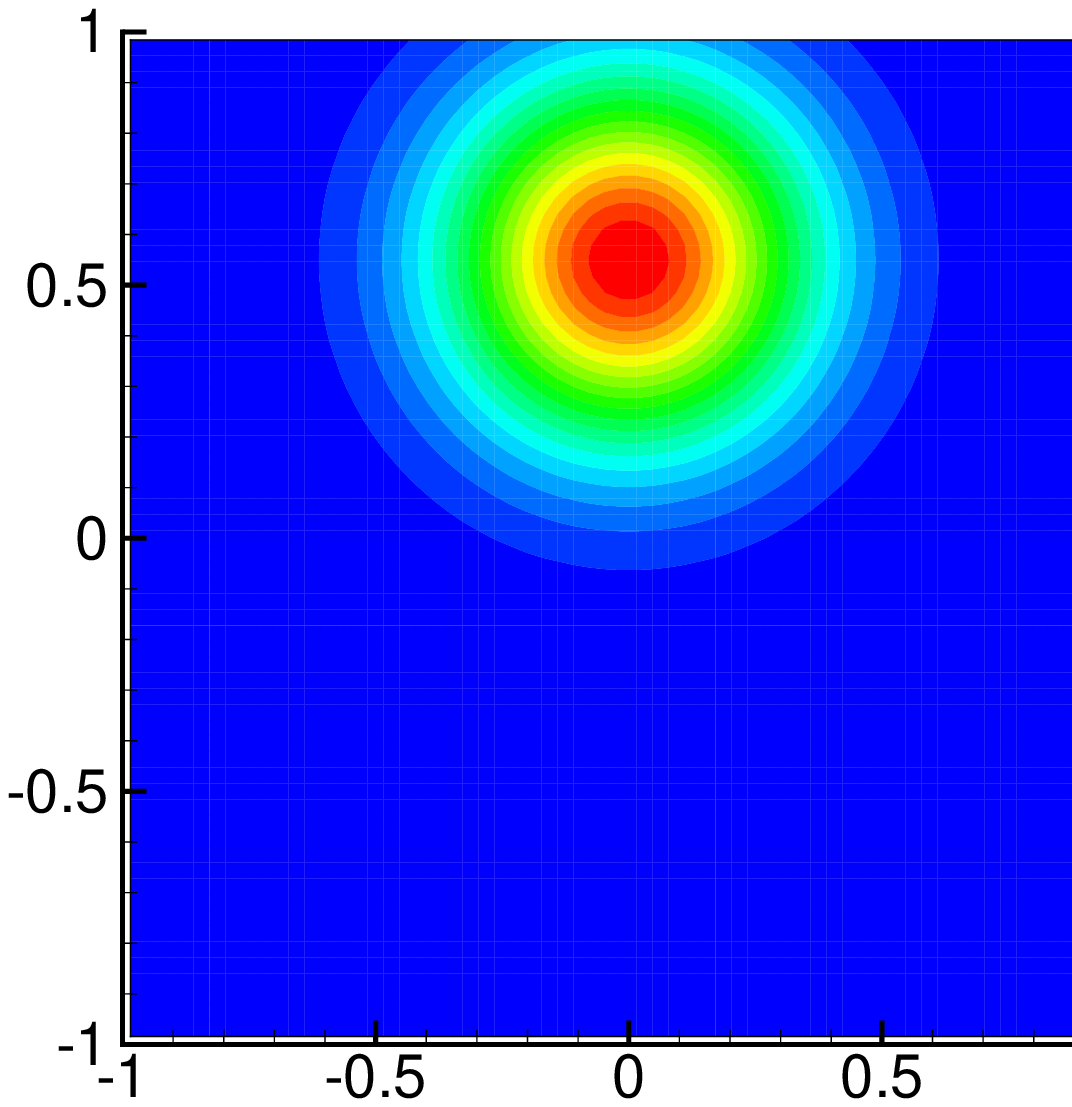}
\caption {Example \ref{ex:2}:  true conductivity.}
\label{fig:true1}
\end{center}
\end{figure}

\begin{figure}
 \begin{center}
    \includegraphics[width = 2in,height=1.8in] {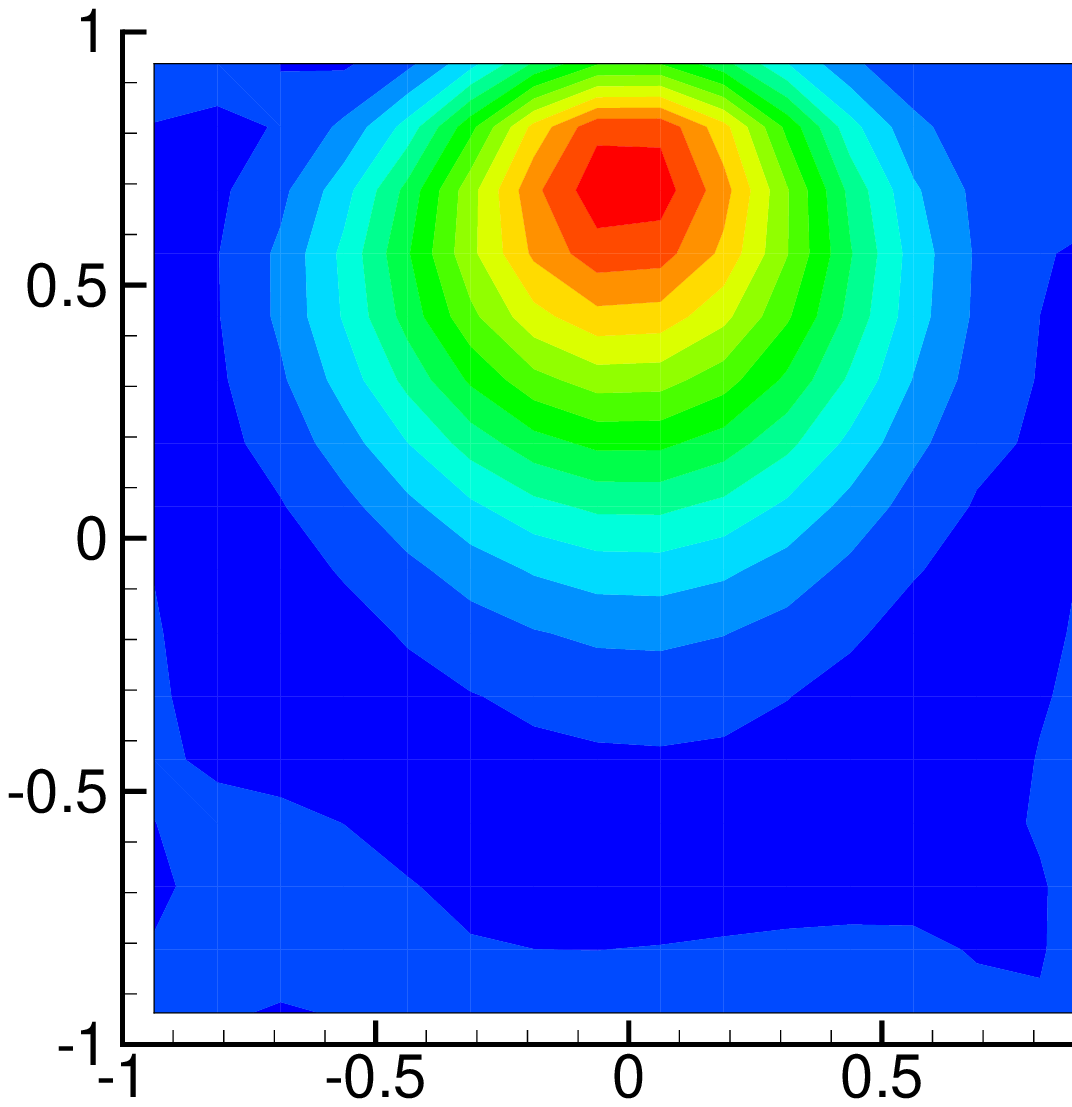}
       \includegraphics[width = 2in,height=1.8in] {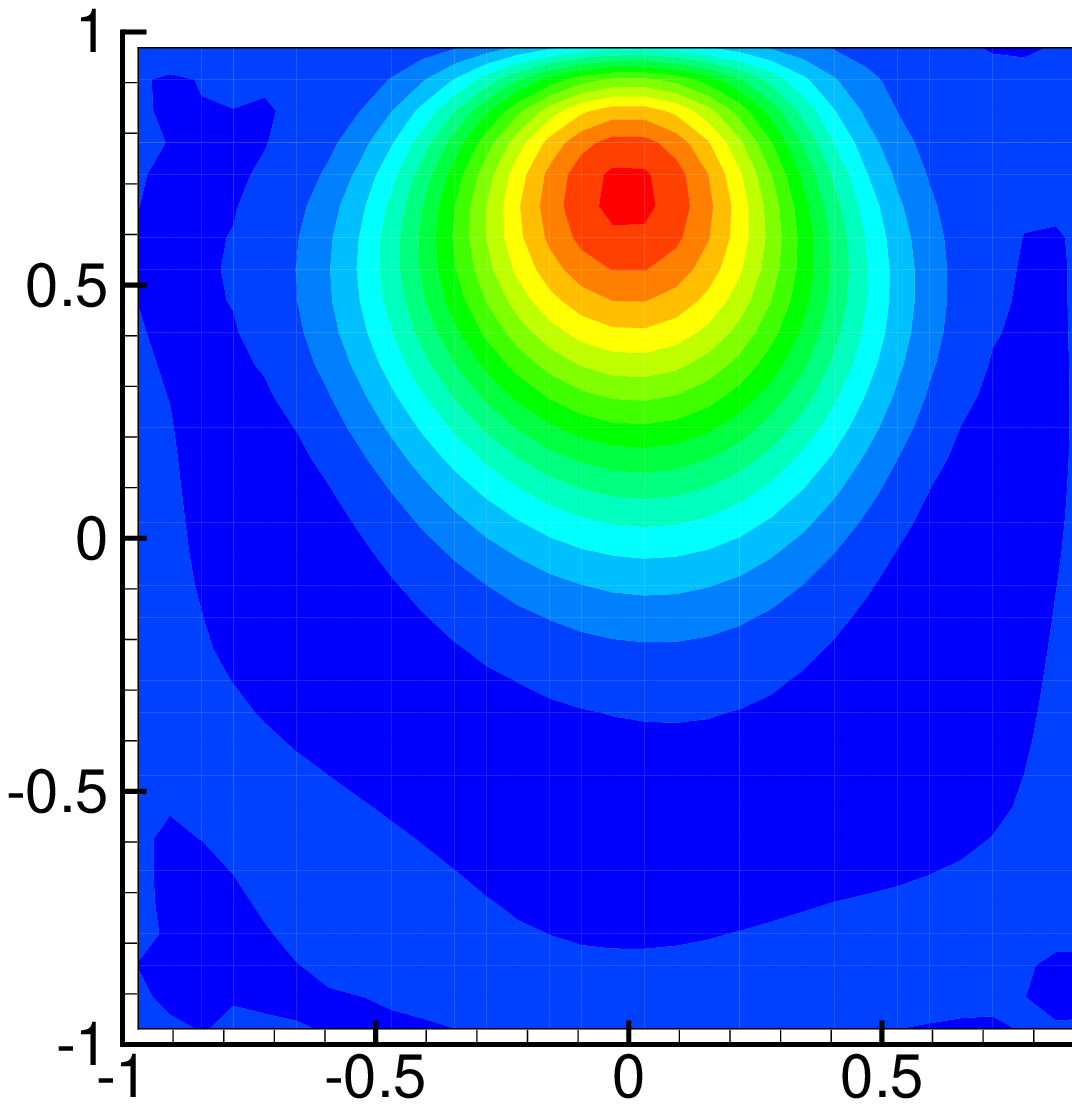}
    \includegraphics[width = 2in,height=1.8in] {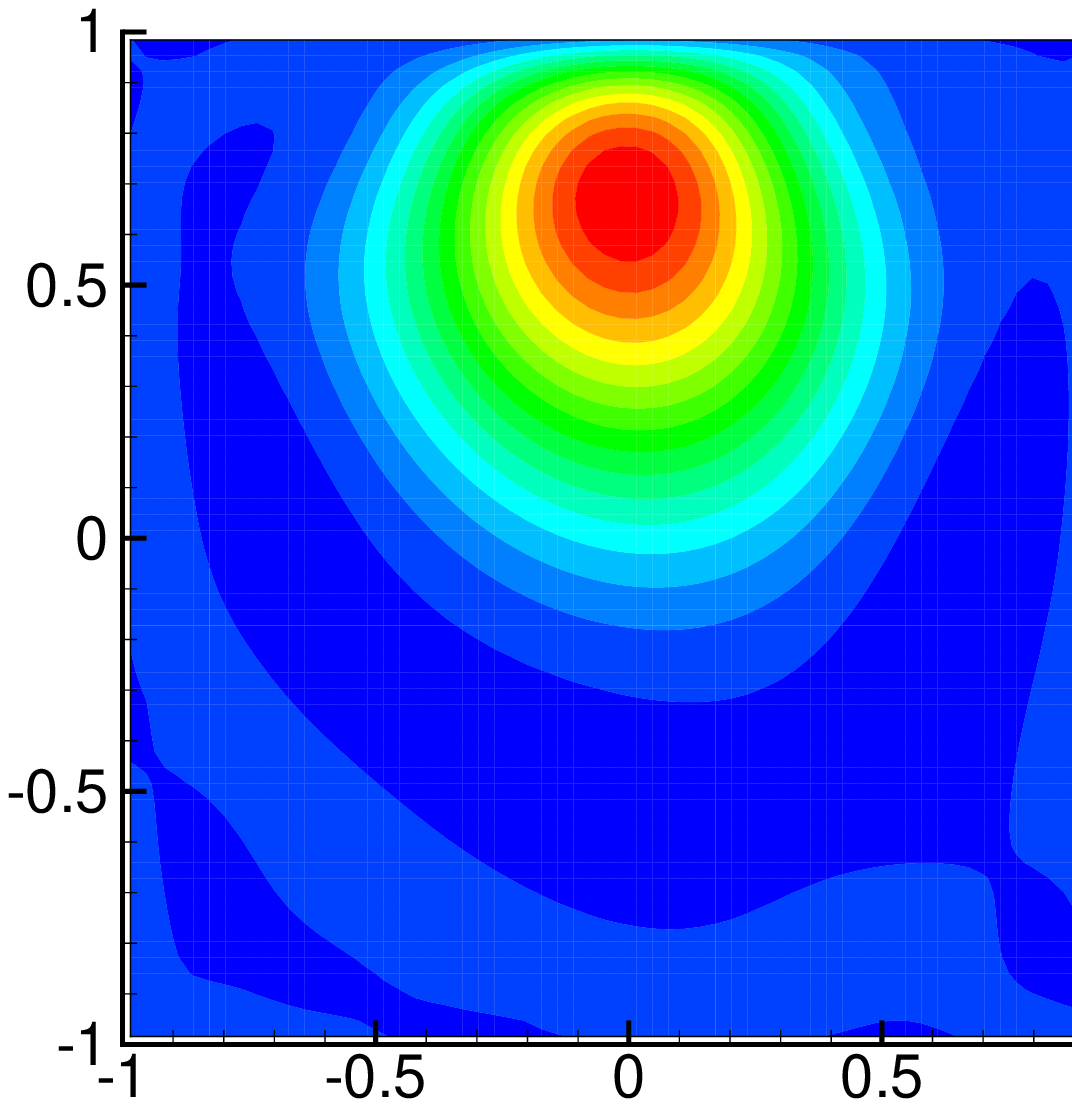}

\caption {Example \ref{ex:2}:  computed conductivity with data noise
$\varepsilon=0\%$. Left: $16\times 16$; Middle: $32\times 32$; Right: $64\times 64$.}
\label{fig:e0_1}
\end{center}
\end{figure}

\begin{figure}
 \begin{center}
    \includegraphics[width = 2in,height=1.8in] {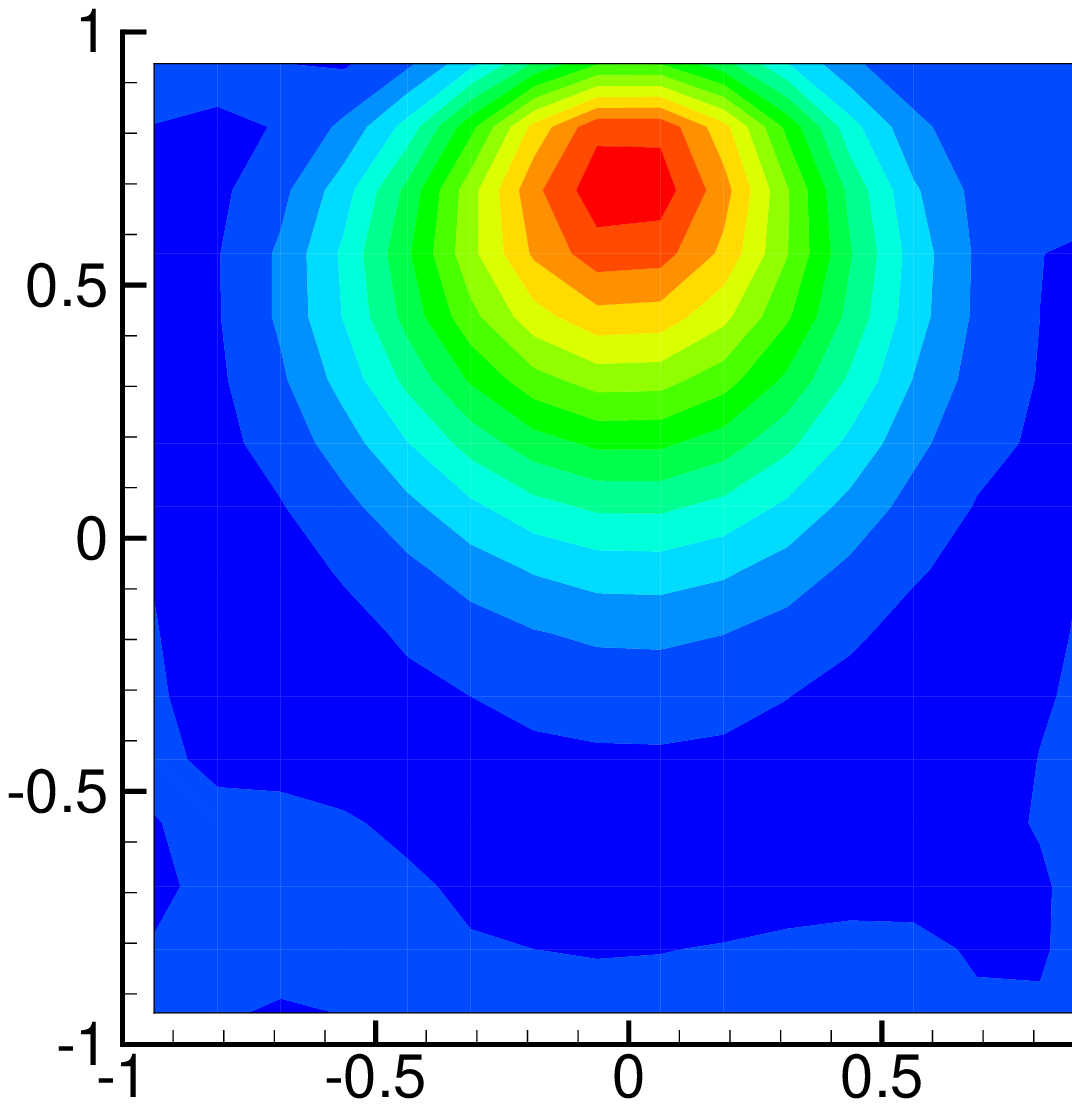}
        \includegraphics[width = 2in,height=1.8in] {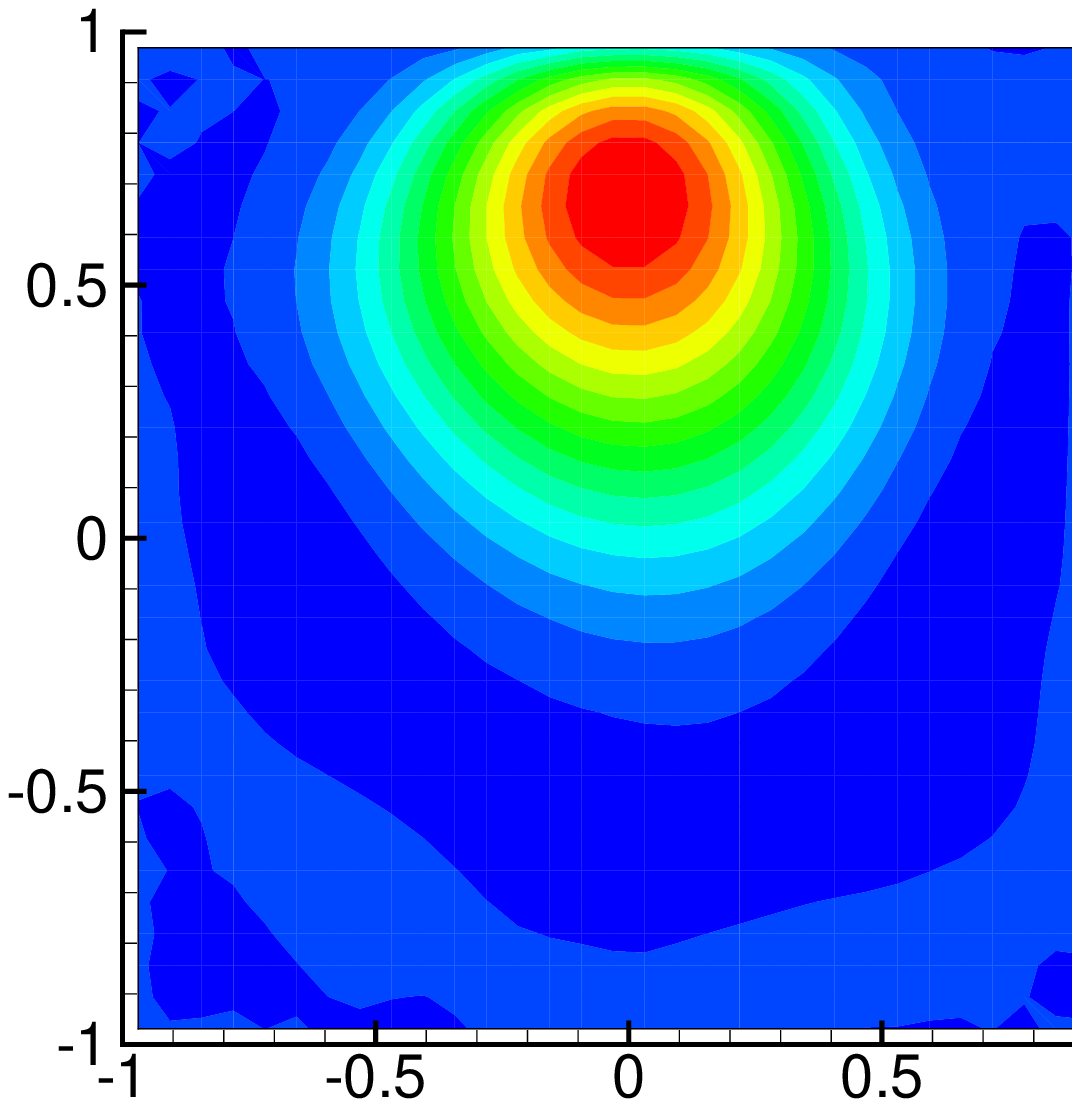}
    \includegraphics[width = 2in,height=1.8in] {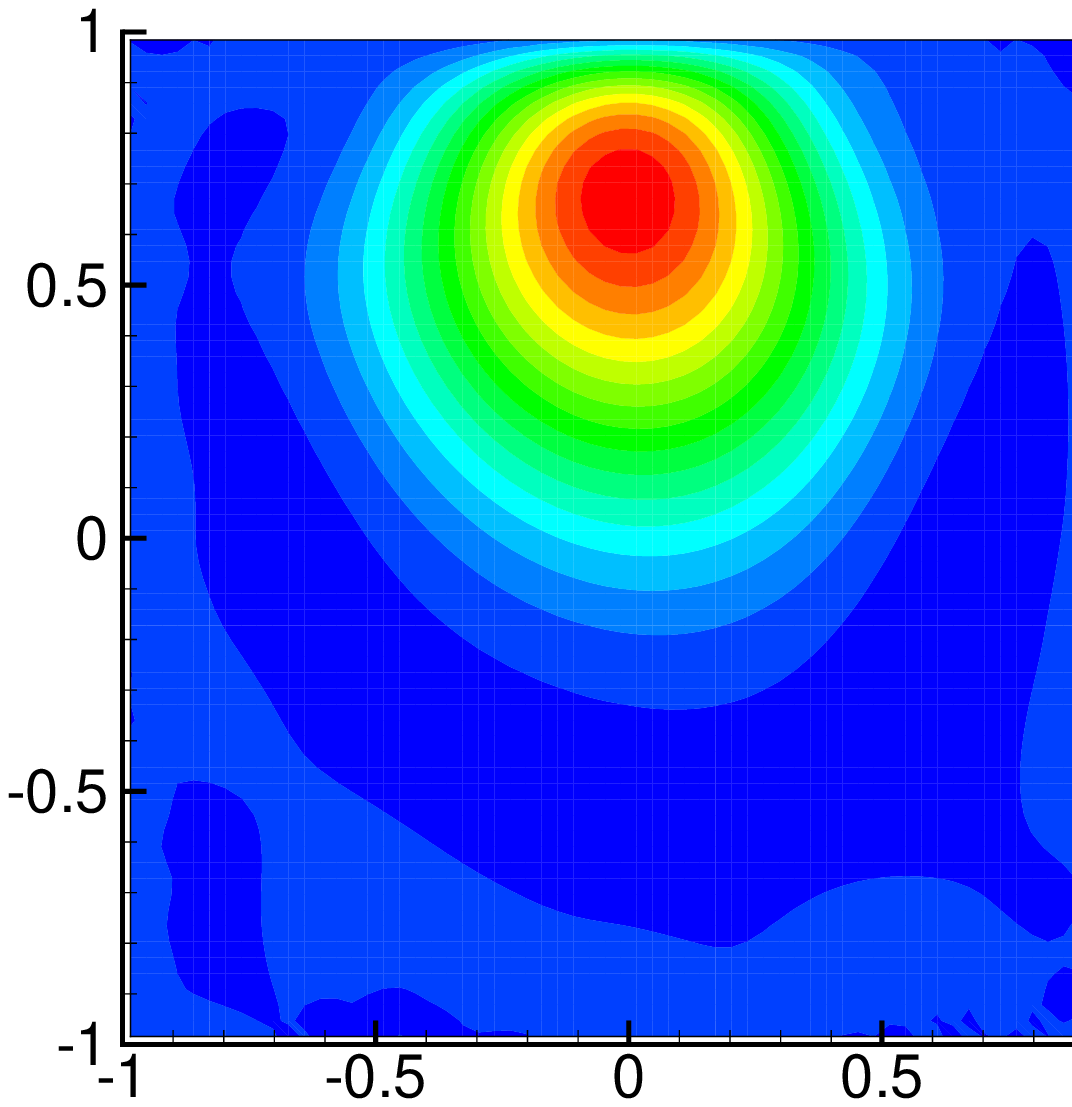}

\caption {Example \ref{ex:2}:  computed conductivity with data noise
$\varepsilon=0.1\%$. Left: $16\times 16$; Middle: $32\times 32$; Right: $64\times 64$.}
\label{fig:e-3_1}
\end{center}
\end{figure}

\begin{figure}
 \begin{center}
    \includegraphics[width = 2in,height=1.8in] {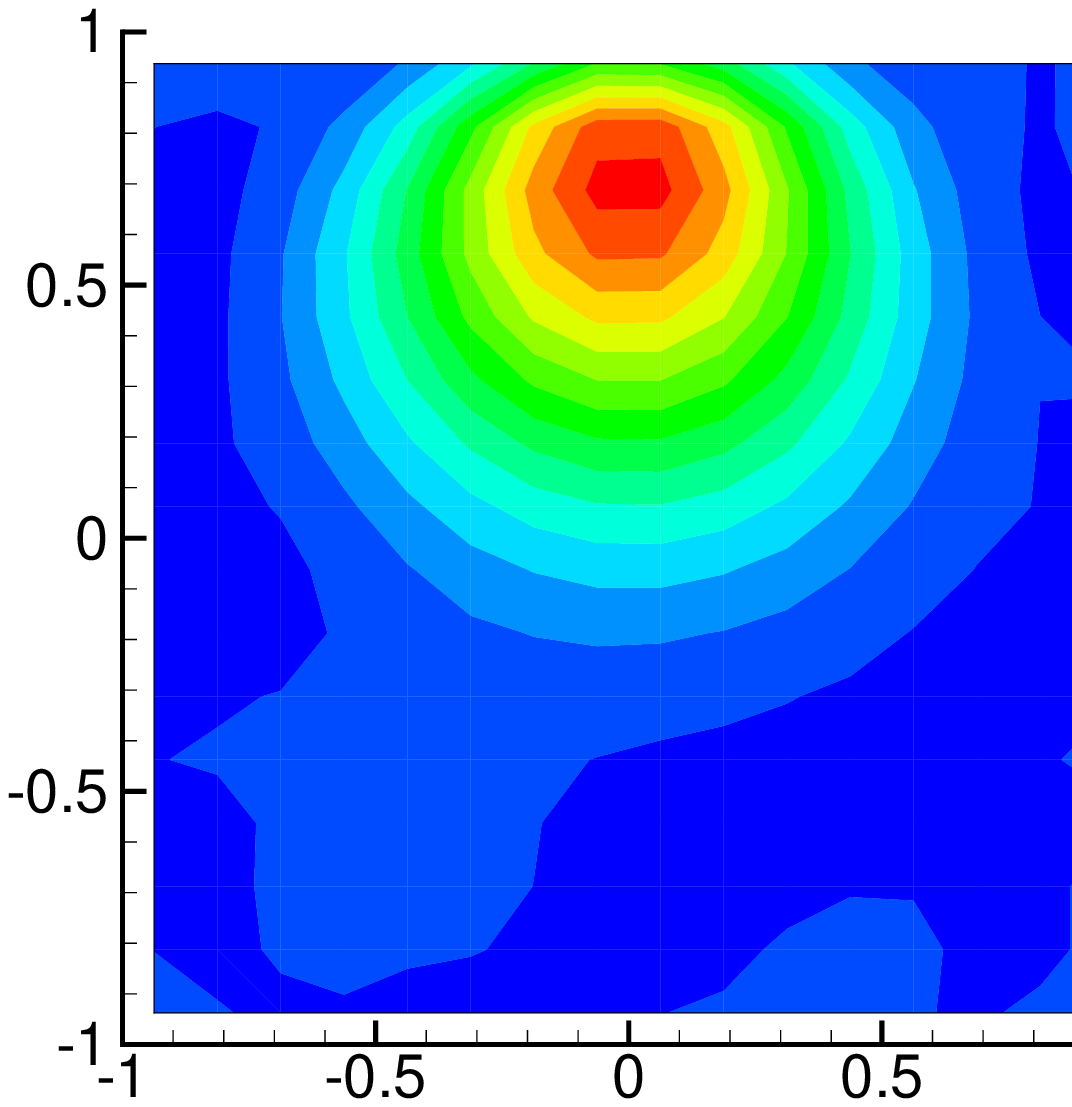}
       \includegraphics[width = 2in,height=1.8in] {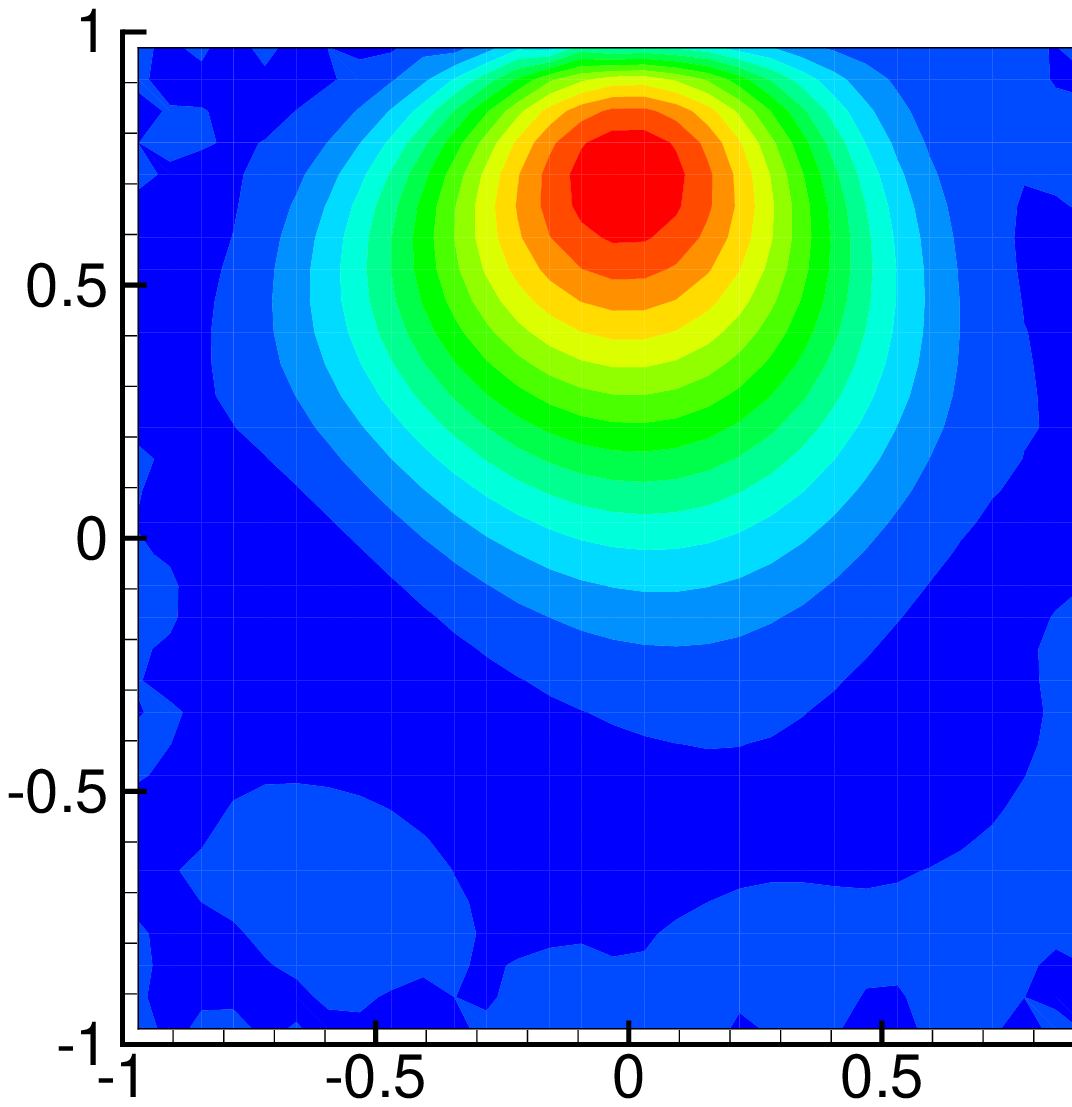}
    \includegraphics[width = 2in,height=1.8in] {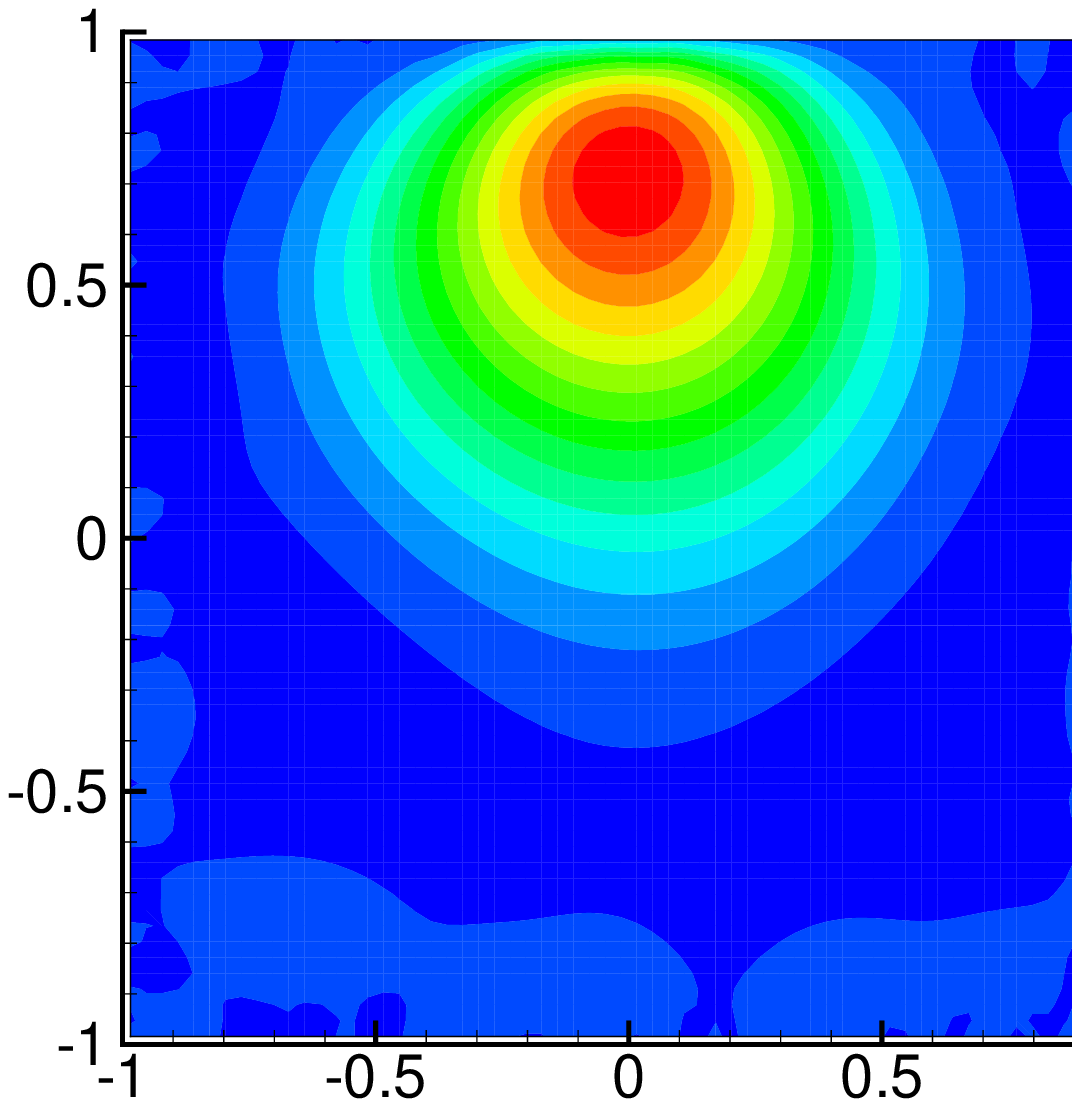}

\caption {Example \ref{ex:2}:  computed conductivity with data noise
$\varepsilon=1\%$. Left: $16\times 16$; Middle: $32\times 32$; Right: $64\times 64$.}
\label{fig:e-2_1}
\end{center}
\end{figure}

\xs{Figure \ref{fig:true1} shows the true conductivity, which has a smooth blob centered at (0, 0.55). We perform our numerical methods by MD-LDG with $P^2$ polynomial space on rectangular meshes. We first set the regularization parameter $\alpha$ to $10^{-8}$, and investigate the effect of various $\alpha$ values afterward. Our study involves three levels of data noise: $\varepsilon=0\%$ (no noise),
$\varepsilon=0.1\%$ and $\varepsilon=1\%$. The computed conductivities under each noise level are presented in Figures  \ref{fig:e0_1}, \ref{fig:e-3_1} and \ref{fig:e-2_1}, respectively. In each group of figures, the mesh sizes are $16\times 16$ (degree of freedom (DOF) 1536), $32\times 32$ (DOF 6144), and $64\times 64$ (DOF 24576) from left to right.
The figures demonstrate that the recoveries from all meshes are able to accurately capture the shape and location of the blob. 
Our result using DOF 6144 is comparable to the adaptive result with DOF 9818 in Example 5.1 of \cite{Zou17}, in terms of similar shape and height of the approximated conductivity (Note that the minimization problem is not the same).
Table  \ref{table:h} lists the heights of the computed conductivities obtained using different meshes and noise levels, where the heights are measured by the maximum value of the conductivity at the centers of all cells.  The true conductivity has a height of 2, and it is apparent that for the same level of noise, finer meshes are able to capture a higher height of the blob and provide a more accurate approximation. Table \ref{table:diff} presents the differences between the computed and measured data $\sum_{j=1}^{M}\|F(\sigma^{computed},f_j)-g_j^\delta\|_{L^2(\partial \Omega)}$. The results indicate that as the mesh is refined, the difference becomes smaller. In general, the results of lower noise levels are better than the results of higher noise levels for the same mesh.
}

\begin{figure}
 \begin{center}
     \includegraphics[width = 2in,height=1.8in] {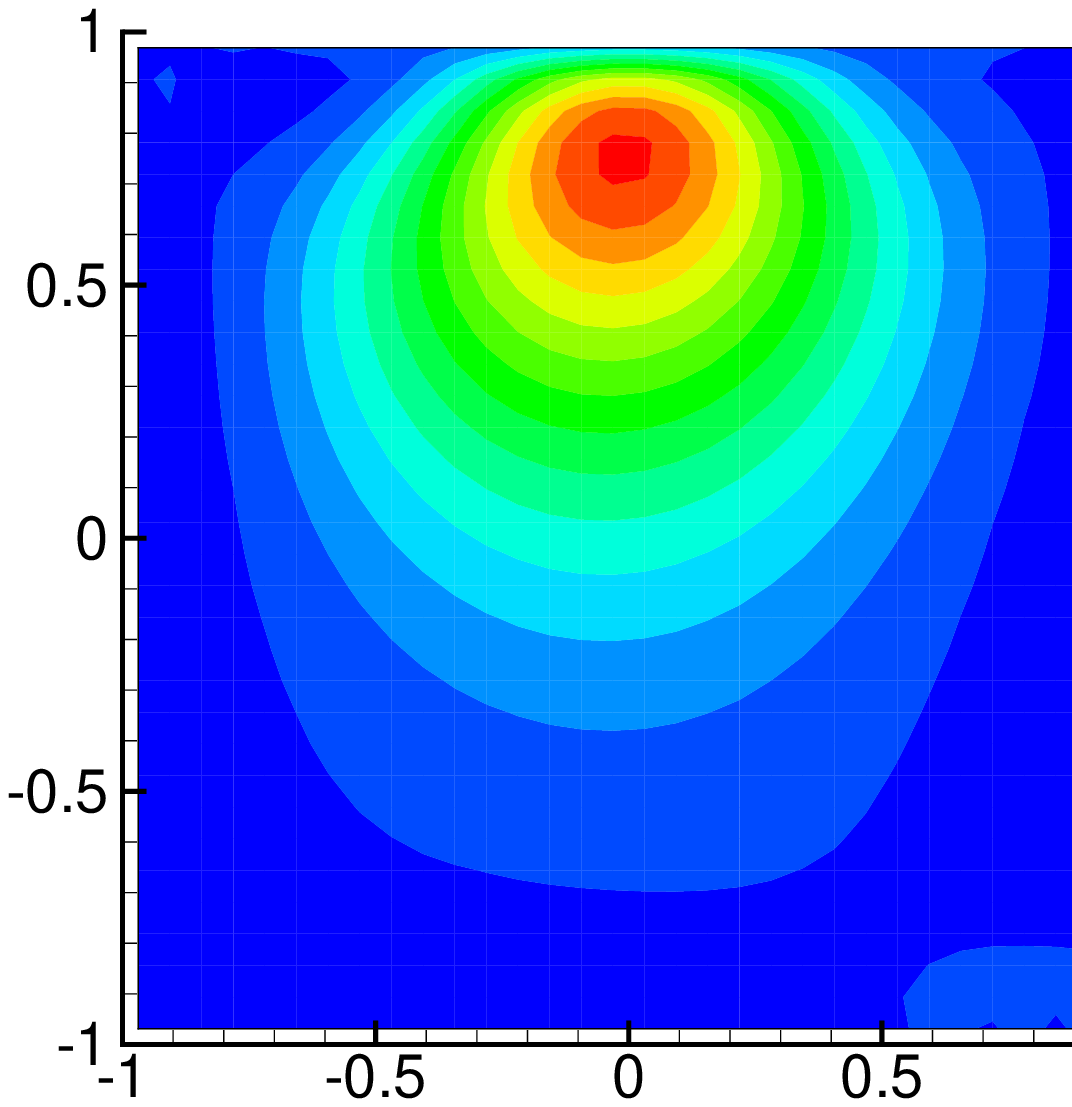}
    \includegraphics[width = 2in,height=1.8in] {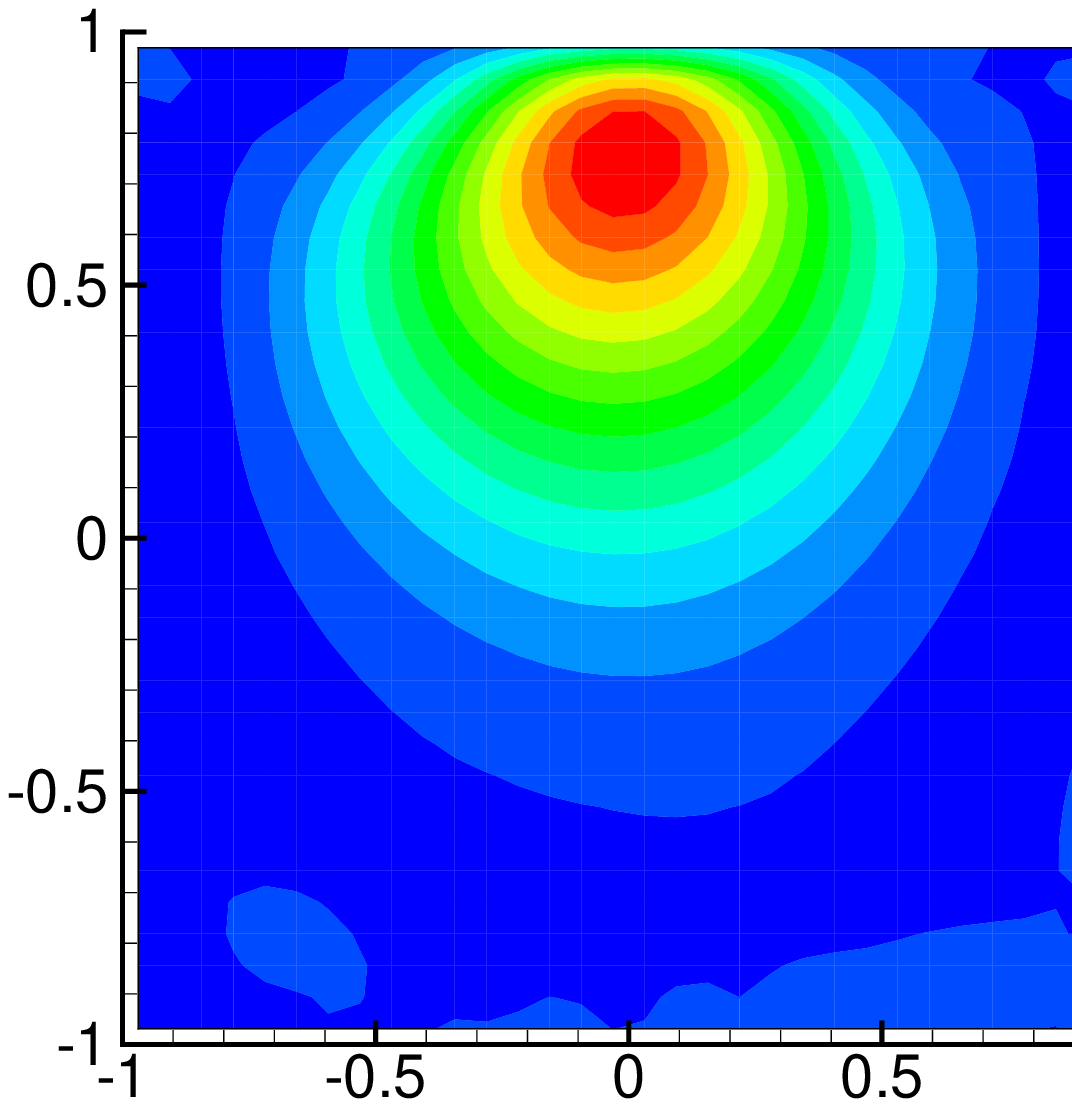}
    \includegraphics[width = 2in,height=1.8in] {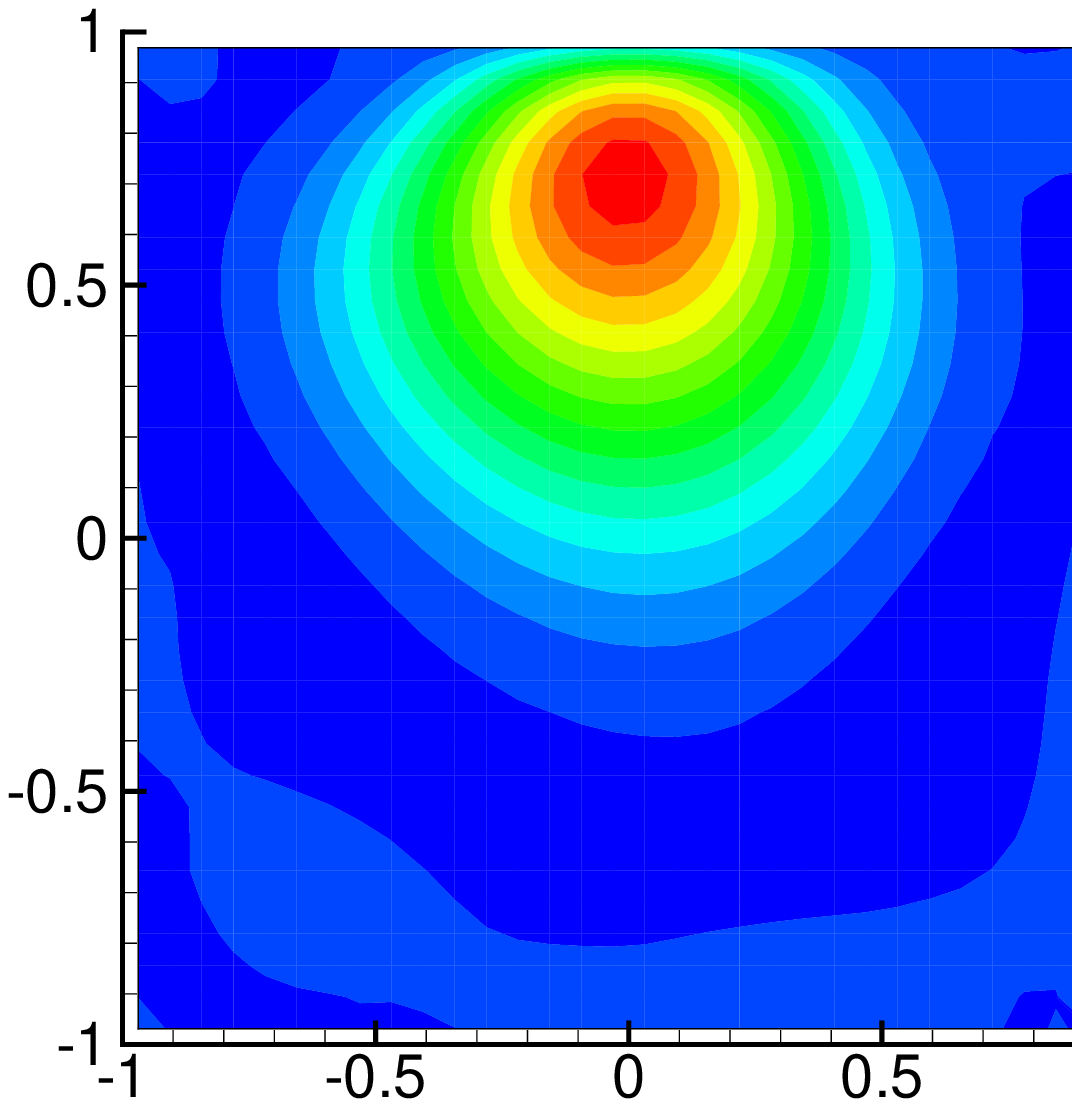}\\
            \includegraphics[width = 2in,height=1.8in] {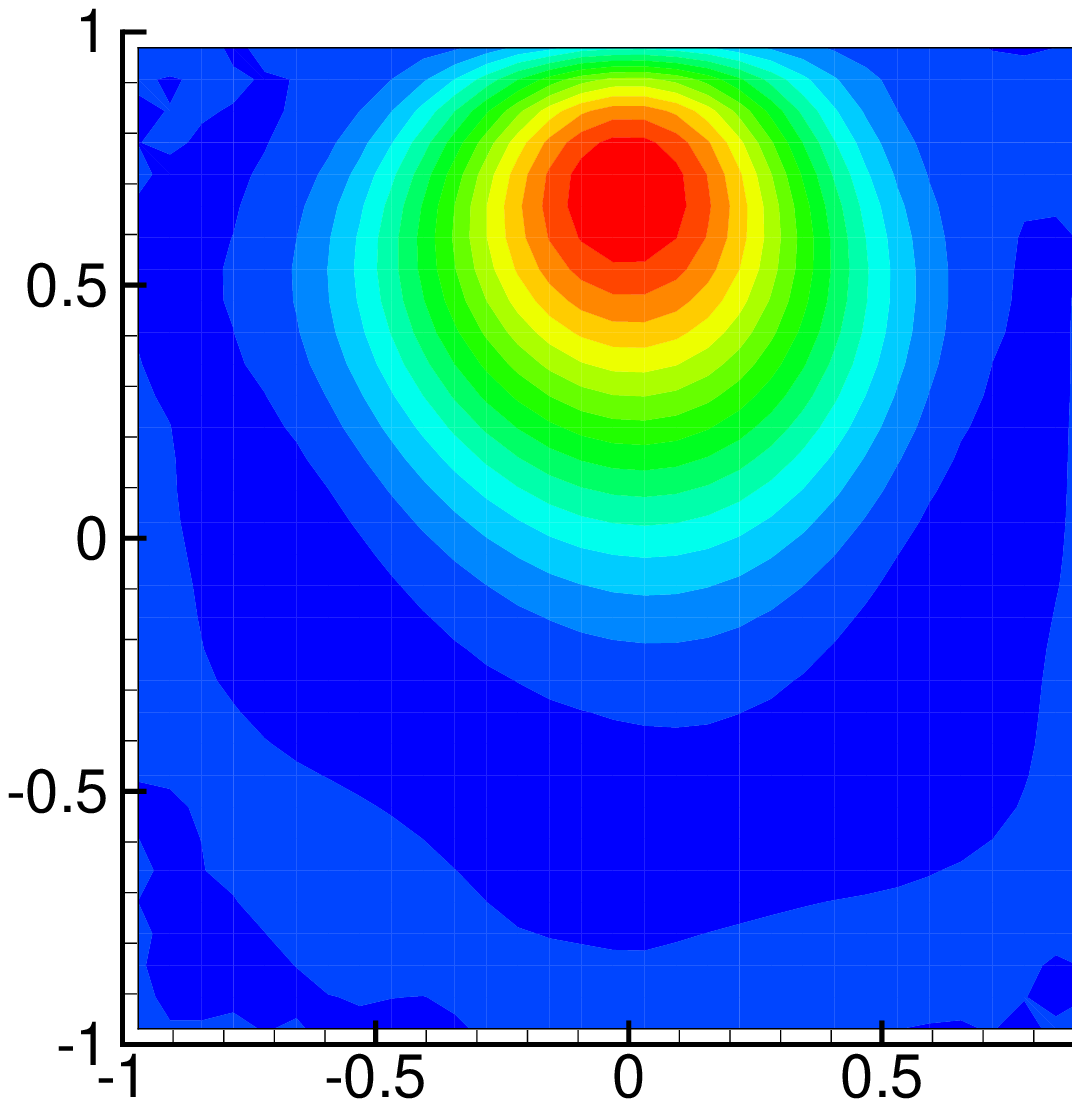}
                \includegraphics[width = 2in,height=1.8in] {n1_a8.eps}
  \includegraphics[width = 2in,height=1.8in] {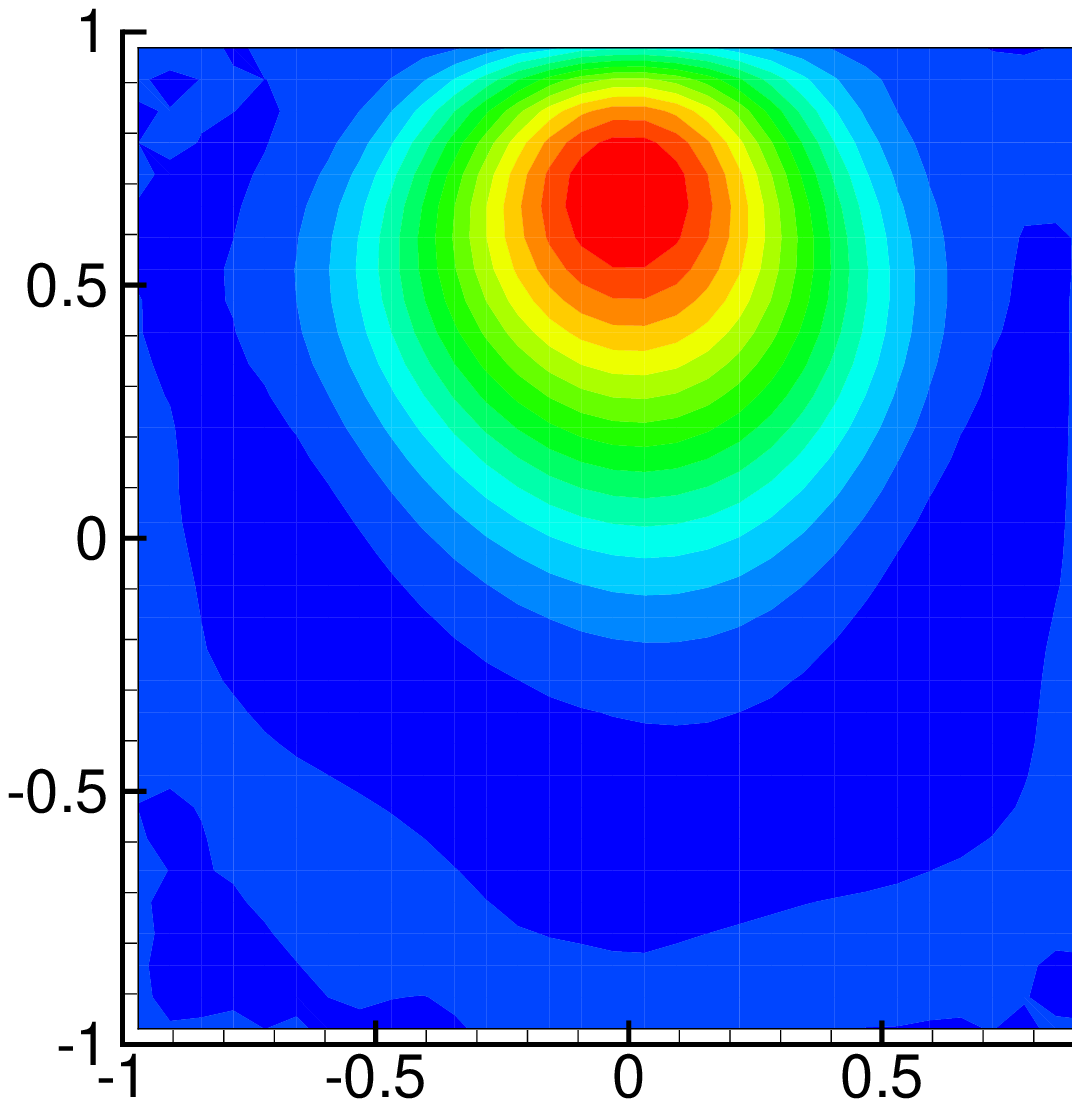}

\caption {Example \ref{ex:2}:  computed conductivity with data noise
$\varepsilon=0.1\%$ and mesh $32\times 32$. Up left: $\alpha=10^{-4}$; Up middle: $\alpha=10^{-5}$, Up right: $\alpha=10^{-6}$;  Down left $\alpha=10^{-7}$; Down middle: $\alpha=10^{-8}$, Down right: $\alpha=0$.}
\label{fig:e-2_a}
\end{center}
\end{figure}

We would like to mention that the proposed method is not sensitive to the regularization parameter $\alpha$.
In Figure \ref{fig:e-2_a}, we show the reconstructions for six different orders of magnitude 
 $\alpha=10^{-4}, 10^{-5}, 10^{-6}, 10^{-7}, 10^{-8}$ and $\alpha=0$. We can see that the reconstructions change
slightly as the regularization parameter varies. Nonetheless, the overall structure of the reconstructions in terms of
conductivity magnitude and center locations remains fairly stable.  From this experiment, we notice that smaller $\alpha$ gives slightly better results with a higher height of the blob and smaller error of conductivity. The results of $\alpha=10^{-8}$ and $\alpha=0$ are indistinguishable. Thus, we will use   $\alpha=10^{-8}$ for all the following examples throughout the paper.
We would like to remark that although we do not see any instability with zero regularization in this particular example, from the analysis we do need a small positive $\alpha$ for stability and convergence.

    \begin{table}
     \centering
   \caption{Example \ref{ex:2}: The  heights of the computed conductivity }
    \smallskip
    \begin{tabular}{cccc}\\\hline\hline
$N$ & $\varepsilon=0\%$ & $\varepsilon=0.1\%$  & $\varepsilon=1\%$   \\\hline\hline

16$\times$16&1.696 &1.681& 1.668      \\
32$\times$32& 1.776&1.758 &  1.695     \\
64$\times$64&1.816 &1.790 & 1.715      \\

\hline\hline
    \end{tabular}\label{table:h}
    \end{table}

    \begin{table}
     \centering
   \caption{Example \ref{ex:2}: The difference between the computed and the measured data $\sum_{j=1}^{M}\|F(\sigma^{computed},f_j)-g_j^\delta\|_{L^2(\partial \Omega)}$ }
    \smallskip
    \begin{tabular}{cccc}\\\hline\hline
$N$ & $\varepsilon=0\%$ & $\varepsilon=0.1\%$  & $\varepsilon=1\%$   \\\hline\hline

16$\times$16&7.74E-2 &7.80E-2& 1.10E-1      \\
32$\times$32& 2.74E-2&2.87E-2 &  8.89E-2     \\
64$\times$64&7.18E-3 &1.10E-2 & 8.08E-2      \\

\hline\hline
    \end{tabular}\label{table:diff}
    \end{table}

\subsection{Example \ref{ex:3}: Reconstruction of EIT: two smooth blobs}
\label{ex:3}

The third example is also a 2D problem  on the domain $\Omega=[-1,1]\times [-1,1]$. 
The true conductivity is given by
$\displaystyle \sigma(x)=\sigma_0(x)+e^{-20((x+0.7)^2+y^2)}+e^{-20(x^2+(y-0.7)^2)}$ with the background
conductivity $\sigma_0=1$, same as  \cite{Zou17}. \xs{We take the background
conductivity as our initial guess.}

\begin{figure}
 \begin{center}
    \includegraphics[width = 2in,height=1.8in] {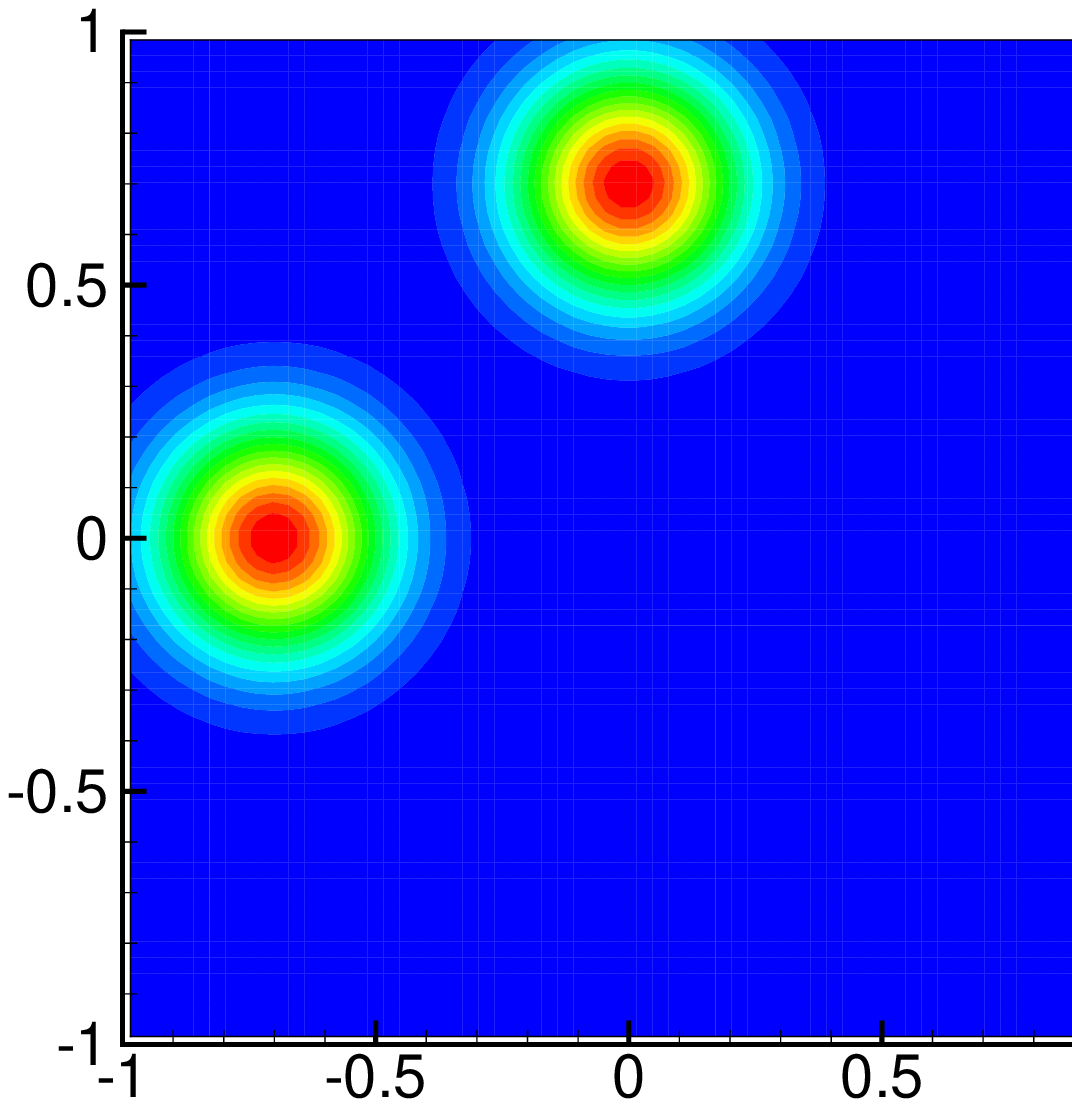}
\caption {Example \ref{ex:3}:  true conductivity.}
\label{fig:true2}
\end{center}
\end{figure}

\begin{figure}
 \begin{center}
    \includegraphics[width = 2in,height=1.8in] {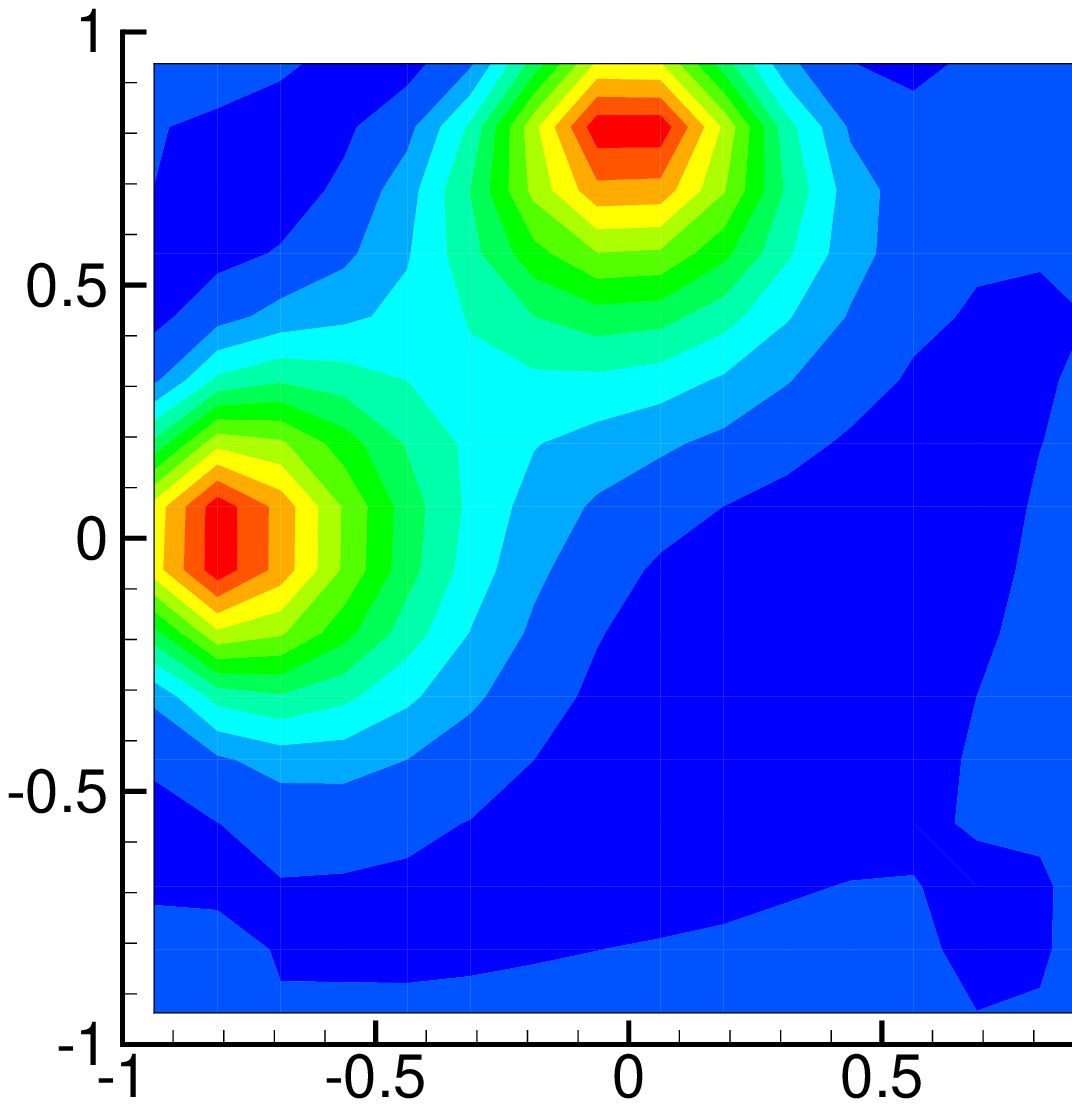}
        \includegraphics[width = 2in,height=1.8in] {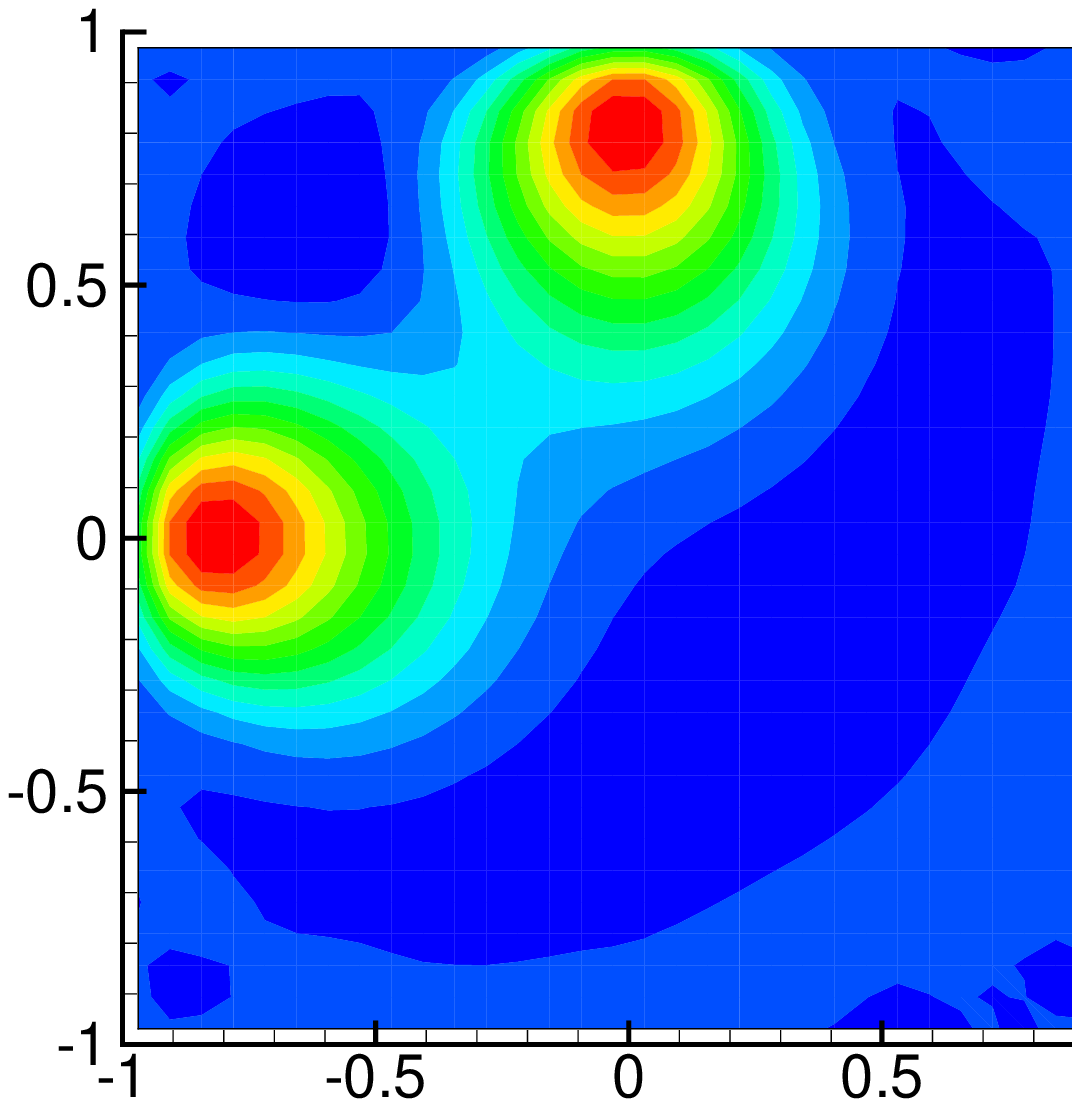}
    \includegraphics[width = 2in,height=1.8in] {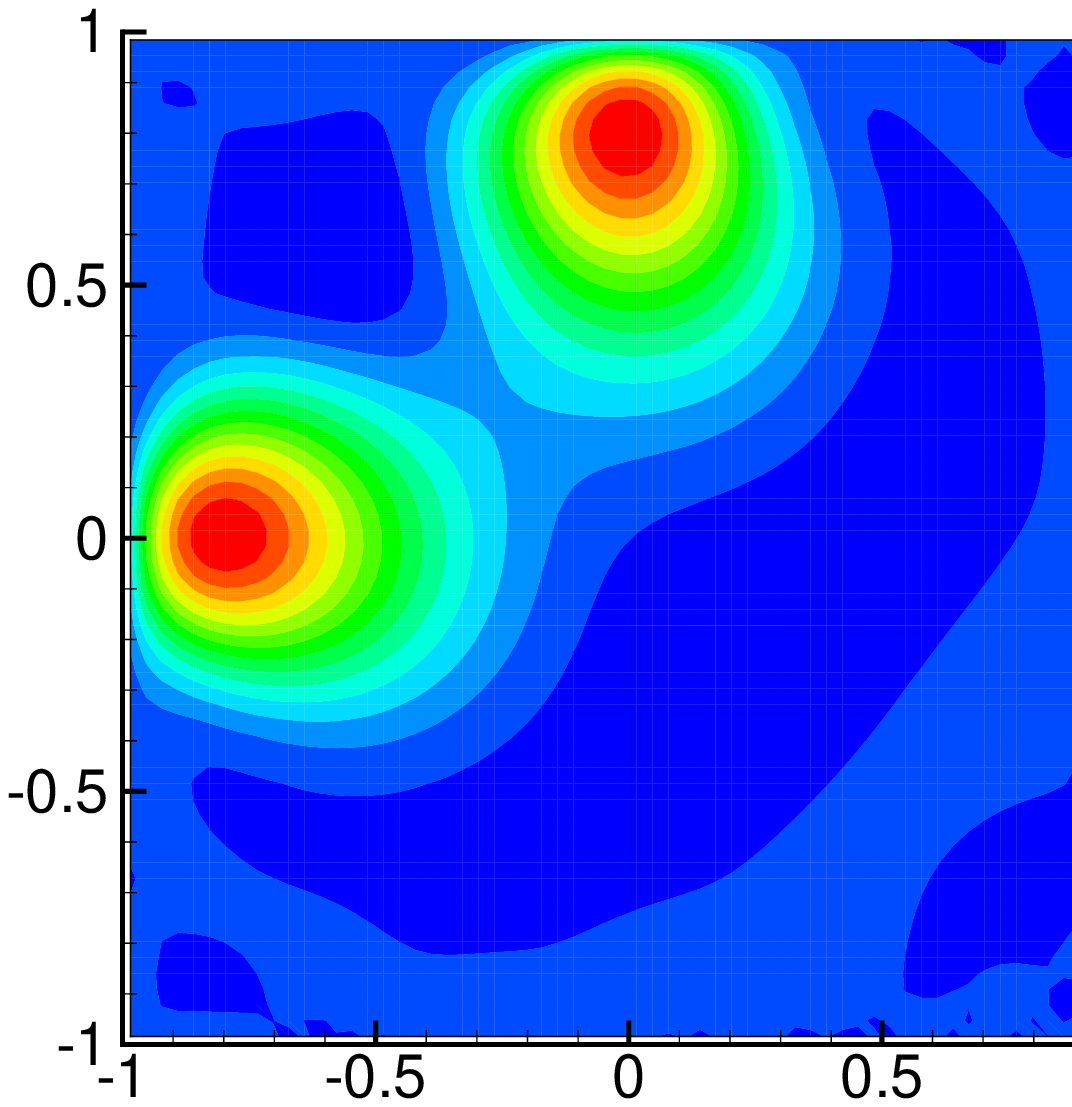}

\caption {Example \ref{ex:3}:  computed conductivity with data noise
$\varepsilon=0.1\%$. Left: $16\times 16$; Middle: $32\times 32$; Right: $64\times 64$.}
\label{fig:e-3_2}
\end{center}
\end{figure}

\begin{figure}
 \begin{center}
    \includegraphics[width = 2in,height=1.8in] {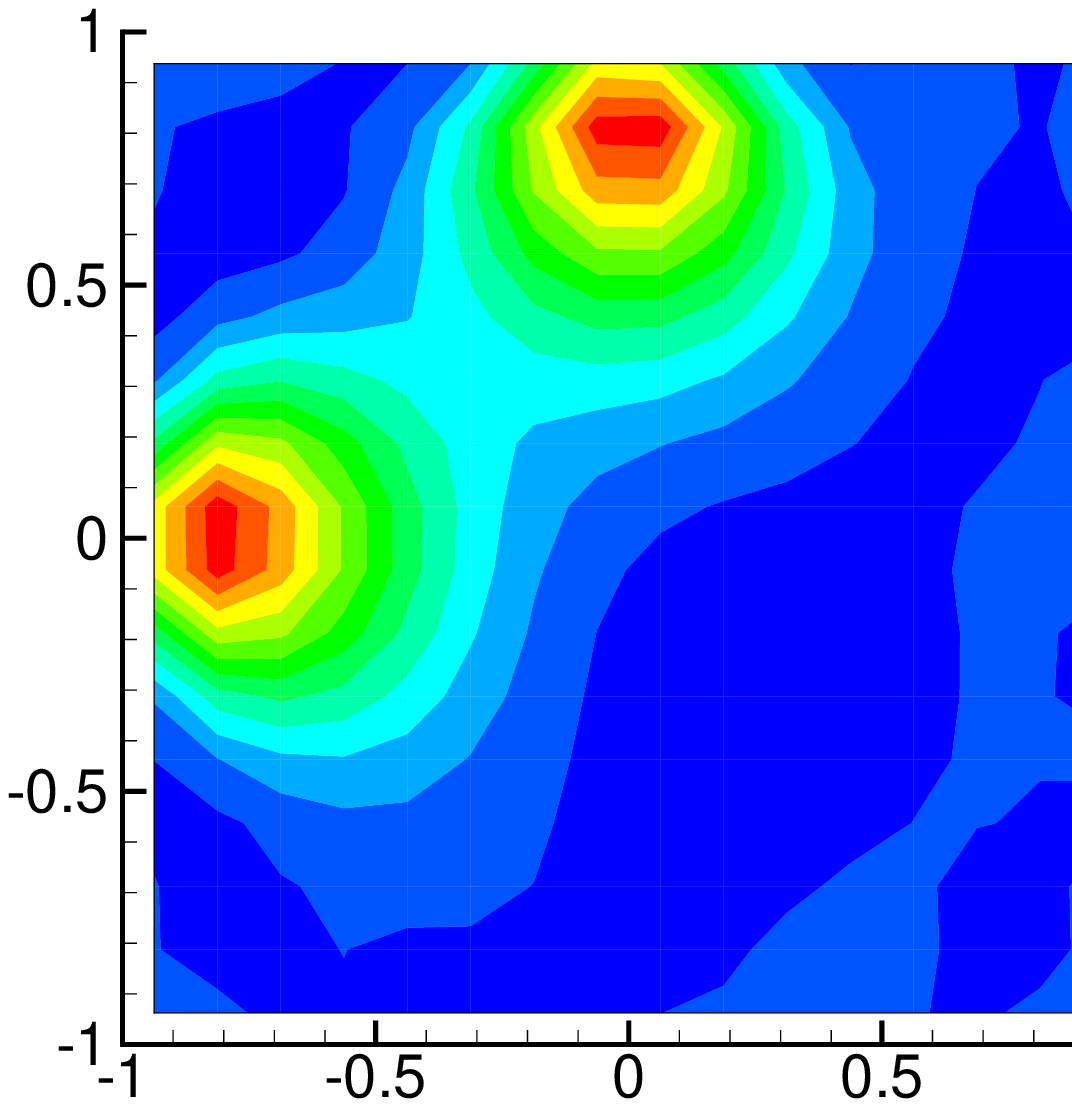}
       \includegraphics[width = 2in,height=1.8in] {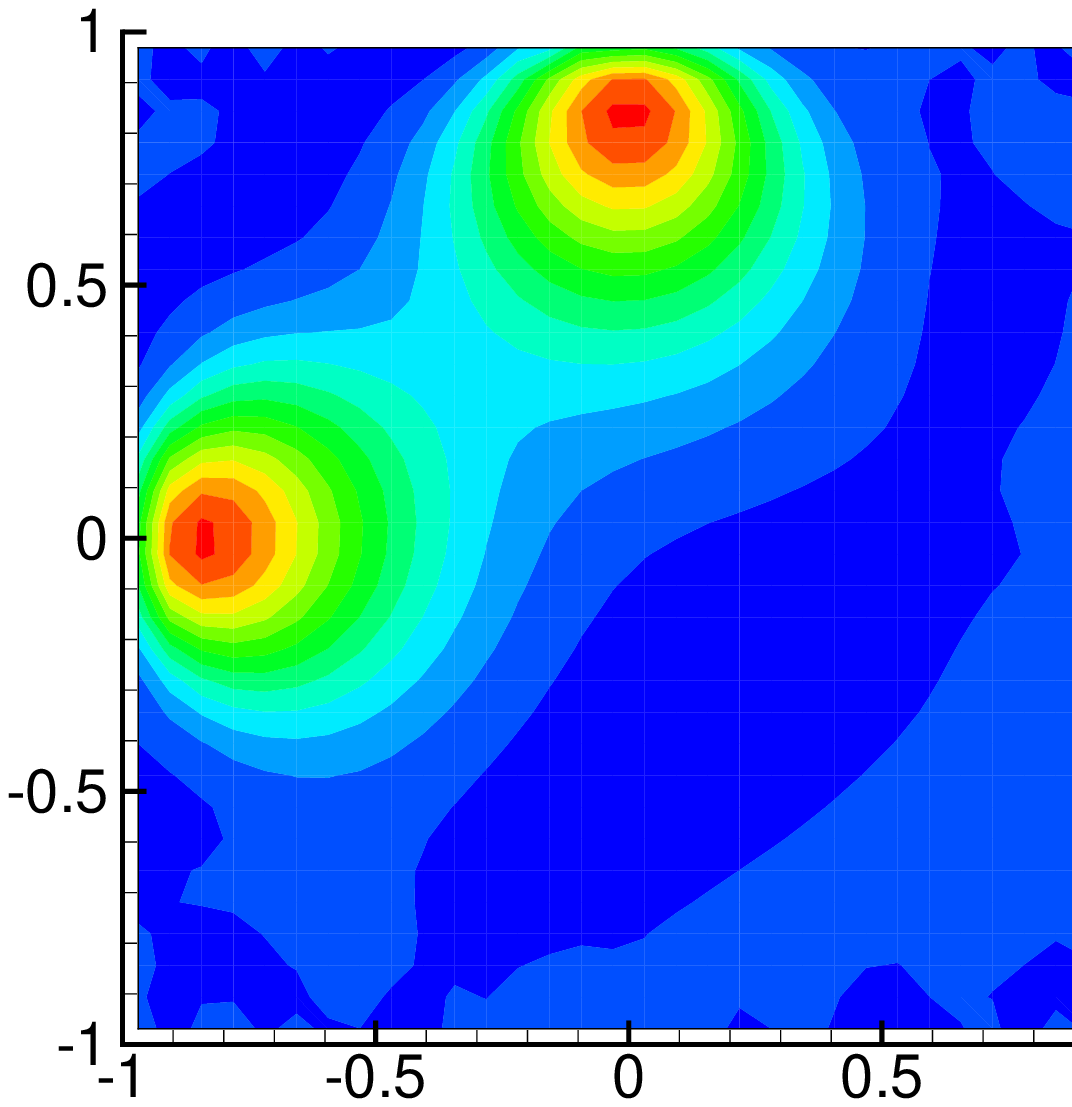}
    \includegraphics[width = 2in,height=1.8in] {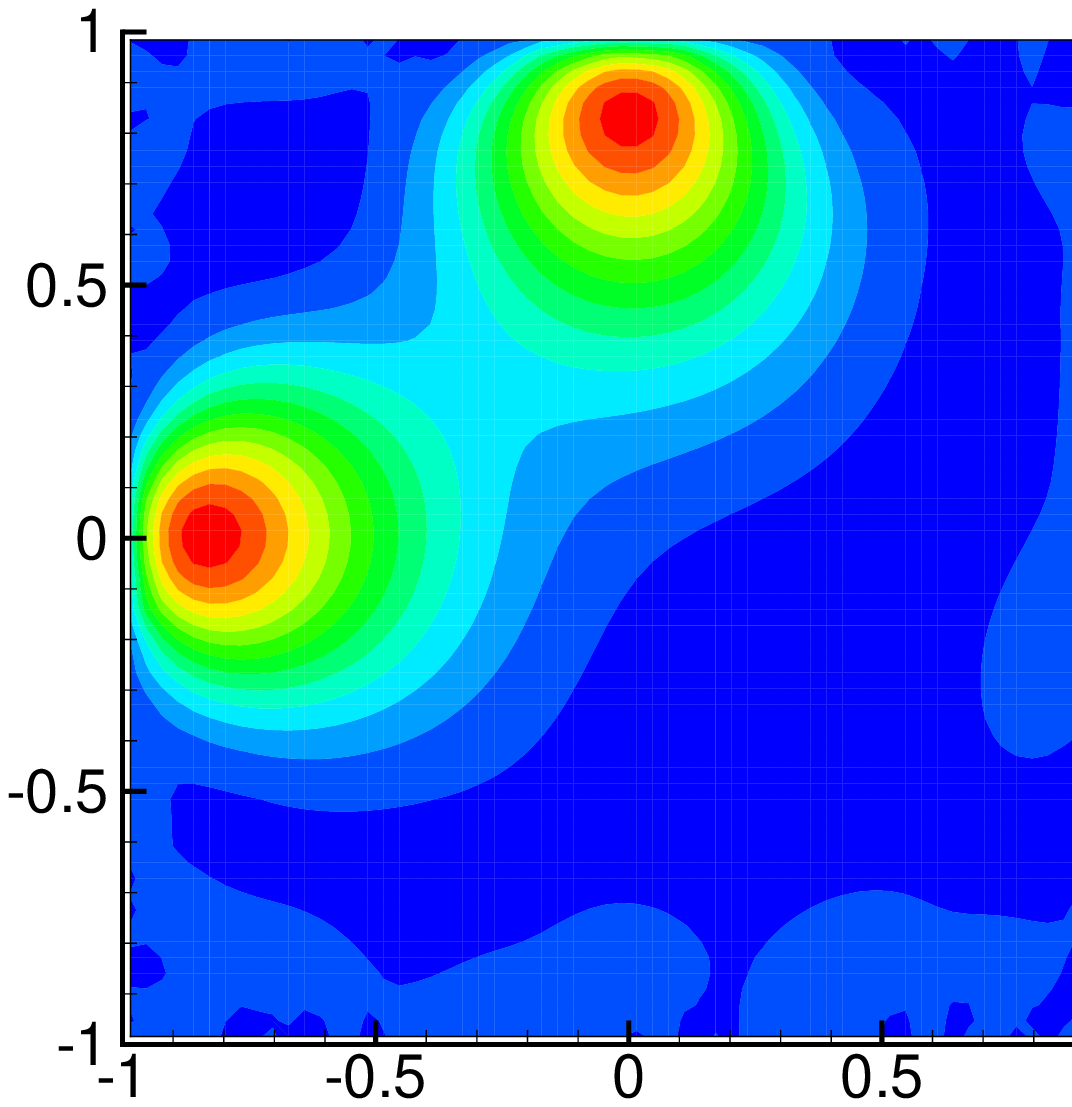}

\caption {Example \ref{ex:3}:  computed conductivity with data noise
$\varepsilon=1\%$. Left: $16\times 16$; Middle: $32\times 32$; Right: $64\times 64$.}
\label{fig:e-2_2}
\end{center}
\end{figure}

The figure of true conductivity is shown in Figure \ref{fig:true2}. It contains two neighboring smooth blobs centered at (-0.7,0) and (0,0.7).
We  consider two levels of data noise
$\varepsilon=0.1\%$ and $\varepsilon=1\%$ and numerical results are computed by MD-LDG with $P^2$ polynomial space on rectangular meshes.
Figures \ref{fig:e-3_2} and \ref{fig:e-2_2} show the computed conductivity with data noise
$\varepsilon=0.1\%$ and $\varepsilon=1\%$ respectively. In both sets of figures, the mesh sizes are $16\times 16$, $32\times 32$, and $64\times 64$  from left to right.
From all the figures, we can see that the recoveries 
capture the location and shape of the two blobs very well. The two blobs are well captured and separated. 
\xs{Our results using fewer DOF are also comparable to the results in Example 5.2 of \cite{Zou17} in terms of similar shape and height of the approximated conductivity.}

\section{{Piecewise continuous conductivity}}

{ In this section, we consider the case when the conductivity is a piecewise continuous function. 
We will redefine the regularity requirement in Section~2 for piecewise continuous conductivity. 
Recall in the minimization functional \eqref{eq_int_R}
\begin{equation*}
R(\sigma)=\frac{1}{2}\sum_{j=1}^{M}\|F(\sigma,f_j)-g_j^\delta\|_{L^2(\partial \Omega)}^2+\frac{\alpha}{2}\|\sigma-\sigma^0\|_{H^1(\Omega)}^2\, ,
\end{equation*}
the first term describes the discrepancy between the measured data and the model-predicted data on the boundary. In Section~2, in order to work on the $L^2(\partial\Omega)$ norm for the easy computation, we {have to} impose some regularity for the conductivity in the definition of the admissible set $\mathcal{A}$. But if the conductivity is a piecewise continuous function, {we can show that} the $L^2(\partial\Omega)$ is well-defined and hence we can release such {a} requirement. }

{More precisely, let $\Omega_m$ $(1\leq m\leq L-1)$ be some pairwise disjoint inclusions in $\Omega$, and denote $\Omega_L=\Omega\setminus \cup_{m=1}^{L-1}\overline\Omega_m$. We suppose that restricting to each 
$\Omega_m$ $(1\leq m\leq L)$, $\sigma(x)\in C^\mu(\Omega_m)$ for some $0<\mu<1$. Clearly, it contains the case that $\sigma$ is a constant on each $\Omega_m$. The following estimate for the conductivity equation \eqref{eq_int_u} are proved in \cite{LV} (see Corollary 1.3 in  \cite{LV})
$$\max_{1\leq m\leq L}\|u\|_{C^{1,\beta}(\overline\Omega_m)}\leq C \|f\|_{C^{1,\beta}(\partial\Omega)}$$
for some $\beta$ $(0<\beta\leq \mu)$, where $C$ may depend on the domain, $\beta$, $C^\mu(\Omega_m)$ norms of $\sigma$, and other factors, but is independent of $f$. We then have
$$\|\nabla u\|_{C^\beta(\partial\Omega)}\leq C \|f\|_{C^{1,\beta}(\partial\Omega)}$$
and
$$\|F(\sigma,f)\|_{L^2(\partial\Omega)}=\|\sigma\frac{\partial u}{\partial \nu}\|_{L^2(\partial\Omega)}\leq C\|\nabla u\|_{C^\beta(\partial\Omega)}\leq C \|f\|_{C^{1,\beta}(\partial\Omega)}.$$
Therefore, the first term of the minimization functional $R(\sigma)$ is well-defined without adding extra smooth conditions on the conductivity. For the regularization term, the $H^1(\Omega)$ norm is used the same as in \eqref{eq_int_R}, which implies that a smoother conductivity is constructed to approximate the true conductivity. The admissible set is now defined as
\begin{equation*}
\mathcal{\tilde{A}}=\{\sigma\in H^{1}(\Omega): 0<c_1<\sigma<c_2, \mbox{ and } \sigma \mbox{ is known on } \partial\Omega\}
\end{equation*}
where $c_1$ and  $c_2$ are fixed numbers. 
Then, the recovery procedure is the same as in Section~4.
}

\subsection{Example \ref{ex:d}: Convergence of forward problem with discontinuous coefficients}
\label{ex:d}
We will first test the convergence of our MD-LDG as the forward solver for model equation \eqref{eq:modeldg}  with discontinuous coefficients.
LDG (including MD-LDG) has the ability to deal with discontinuous coefficients as long as the mesh is aligned with
the discontinuous interface.

We take the example from \cite{Yanjue}. The computational domain is a square $[0,1]\times[0,1]$.
The coefficient $\sigma$ is a piecewise constant
\begin{equation}
\sigma=\left\{
\begin{array}{ll}
1,& x<0.5 \\
10, & x>0.5.
\end{array}\right.
\end{equation}
We choose the exact solution to be
\[u=\frac{1}{\sigma}\sin\left(\frac{\pi x}{2}\right)(x-0.5)(y-0.5)(x^2+y^2+1).\]
The right hand side $r(x,y)$ and the boundary $b(x,y)$ in \eqref{eq:modeldg} are provided from the calculation of $u$.
We use the MD-LDG with $P^2$ polynomial space.
Table \ref{table:2} showed the $L^2$-errors and orders of accuracy of 
$u$, $\sigma u_x$ and $\sigma u_y$. We again see third order convergence for $u$ and second order for $\sigma u_x$ and $\sigma u_y$.
      \begin{table}
     \centering
   \caption{Example \ref{ex:d}: $L^2$-errors and orders of accuracy of MD-LDG $P^2$. }
    \smallskip
    \begin{tabular}{ccccccc}\\\hline\hline
&$u$ && $\sigma u_x$& & $\sigma u_y$&\\\hline
$N$ & error & order  & error & order & error & order \\\hline\hline
8$\times$ 8&1.08E-04   &--  & 2.94E-03   &--  & 2.68E-03& --   \\

16$\times$16& 1.28E-05 &3.08    &7.33E-04   &  2.00& 6.89E-04& 1.96   \\
32$\times$32& 1.54E-06 & 3.06   & 1.81E-04  & 2.02 & 1.75E-04&1.98    \\
64$\times$64&1.89E-07  &3.02    &  4.50E-05 & 2.01 &4.41E-05 &1.99    \\

\hline\hline
    \end{tabular}\label{table:2}
    \end{table}

\subsection{{Example \ref{ex:d2}: Reconstruction of EIT: inclusions with a constant background}}
\label{ex:d2}
\begin{figure}
 \begin{center}
    \includegraphics[width = 2.in,height=1.7in] {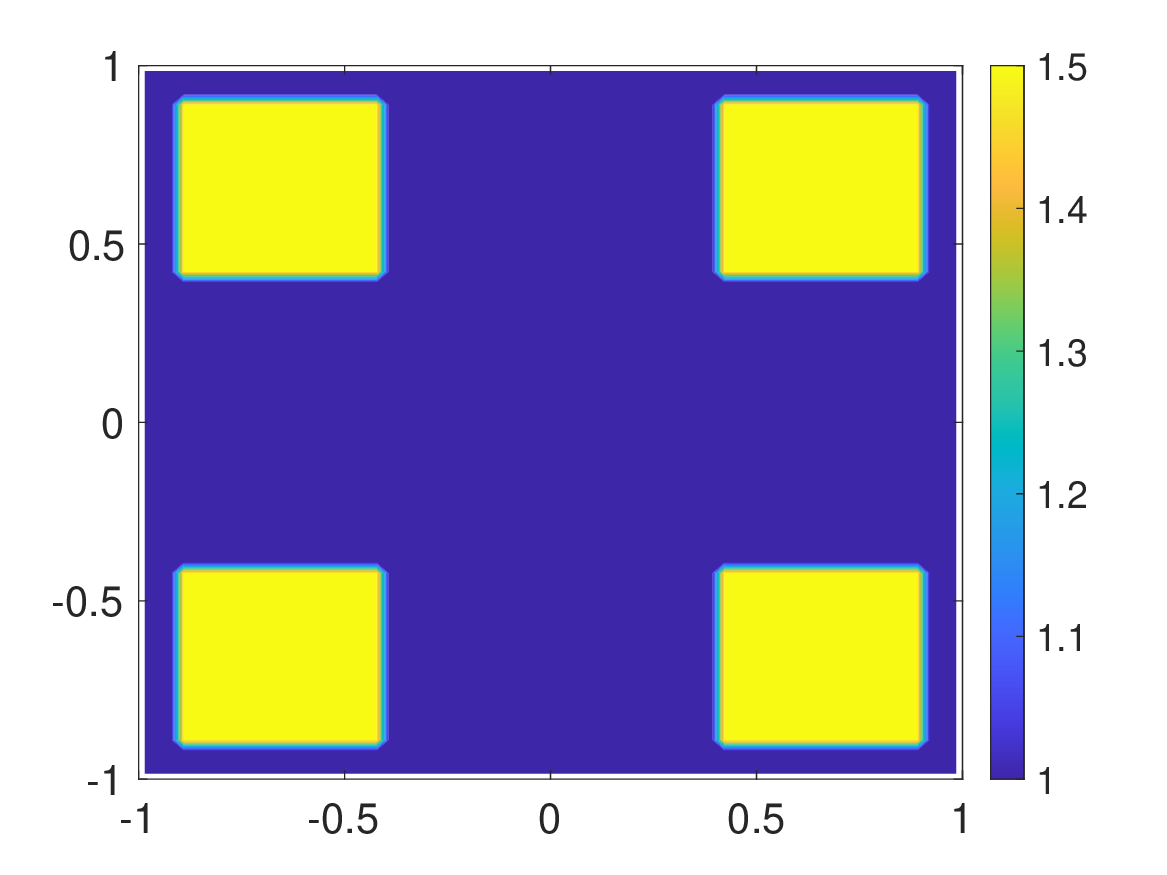}
\caption {Example \ref{ex:d2}:  true conductivity.}
\label{fig:true7}
\end{center}
\end{figure}

\begin{figure}
 \begin{center}
    \includegraphics[width = 2.in,height=1.7in] {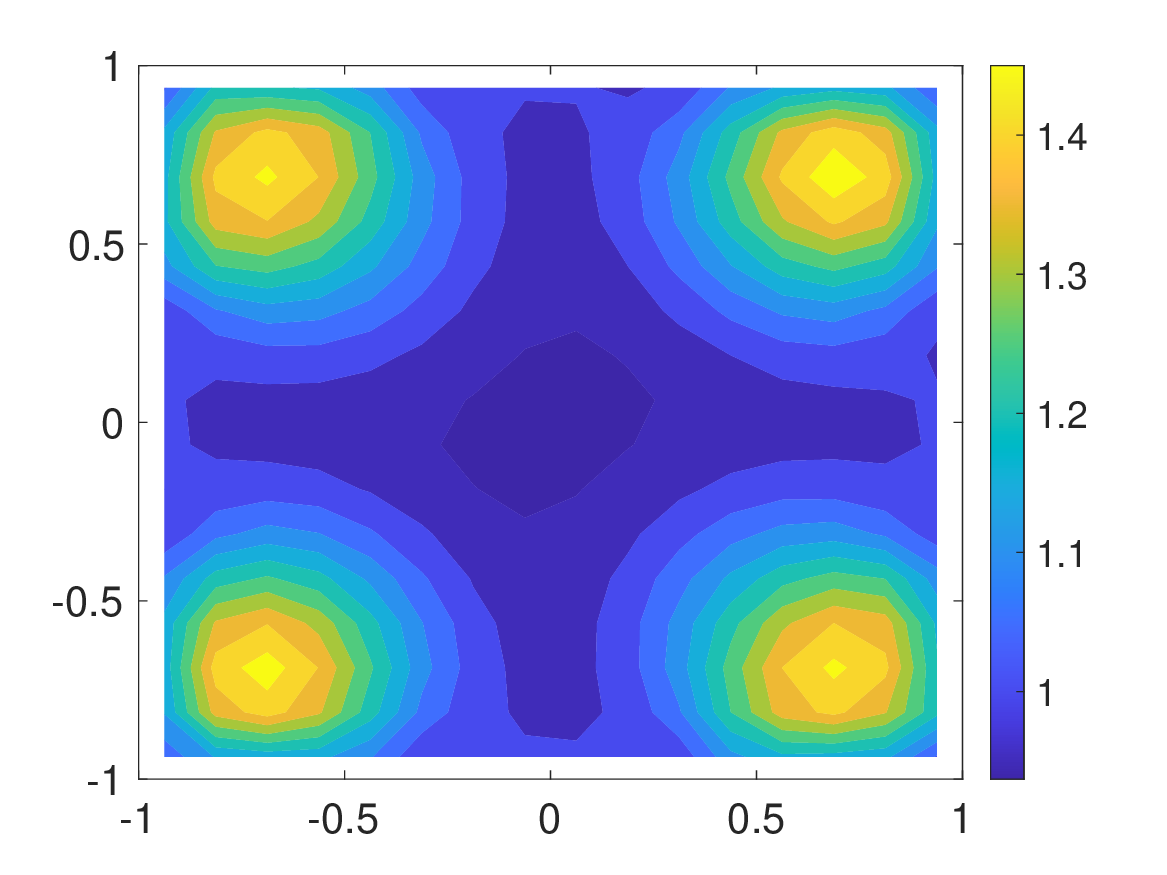}
        \includegraphics[width = 2.in,height=1.7in] {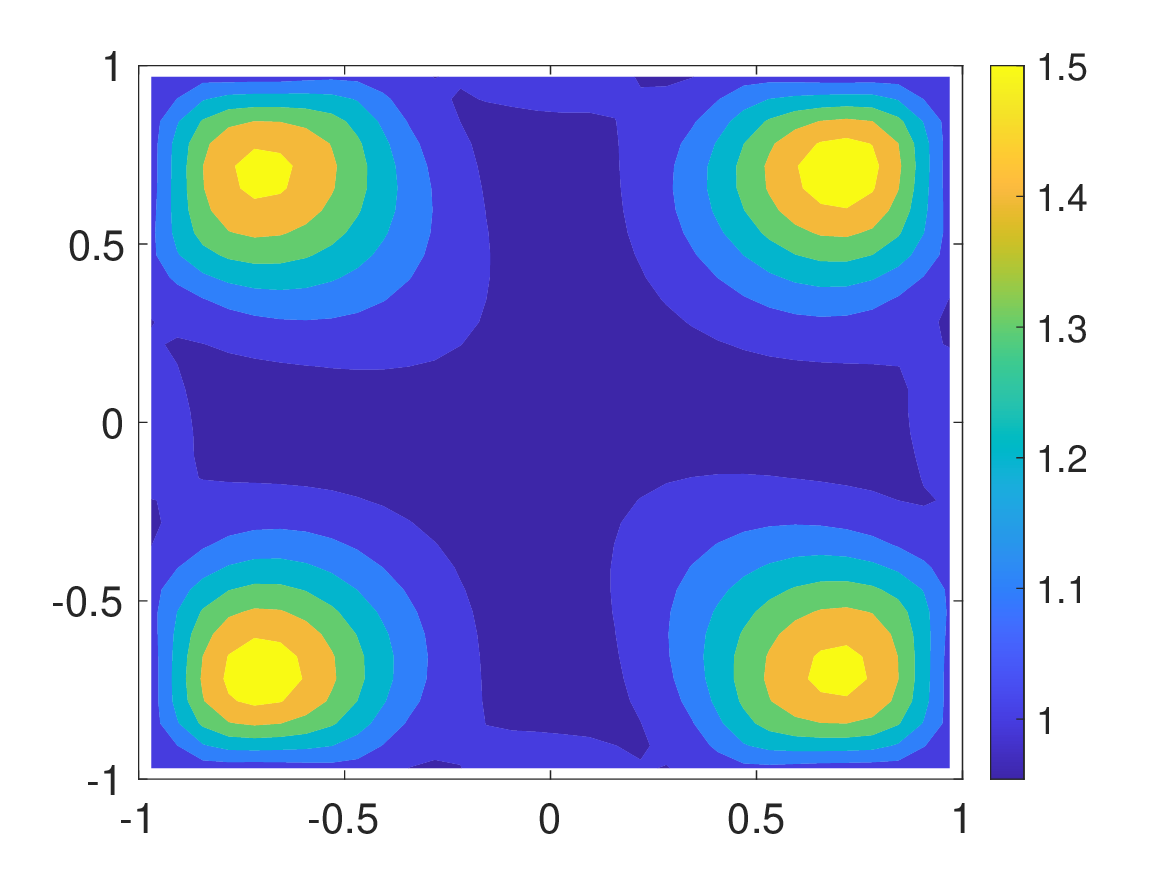}
   \includegraphics[width = 2.in,height=1.7in] {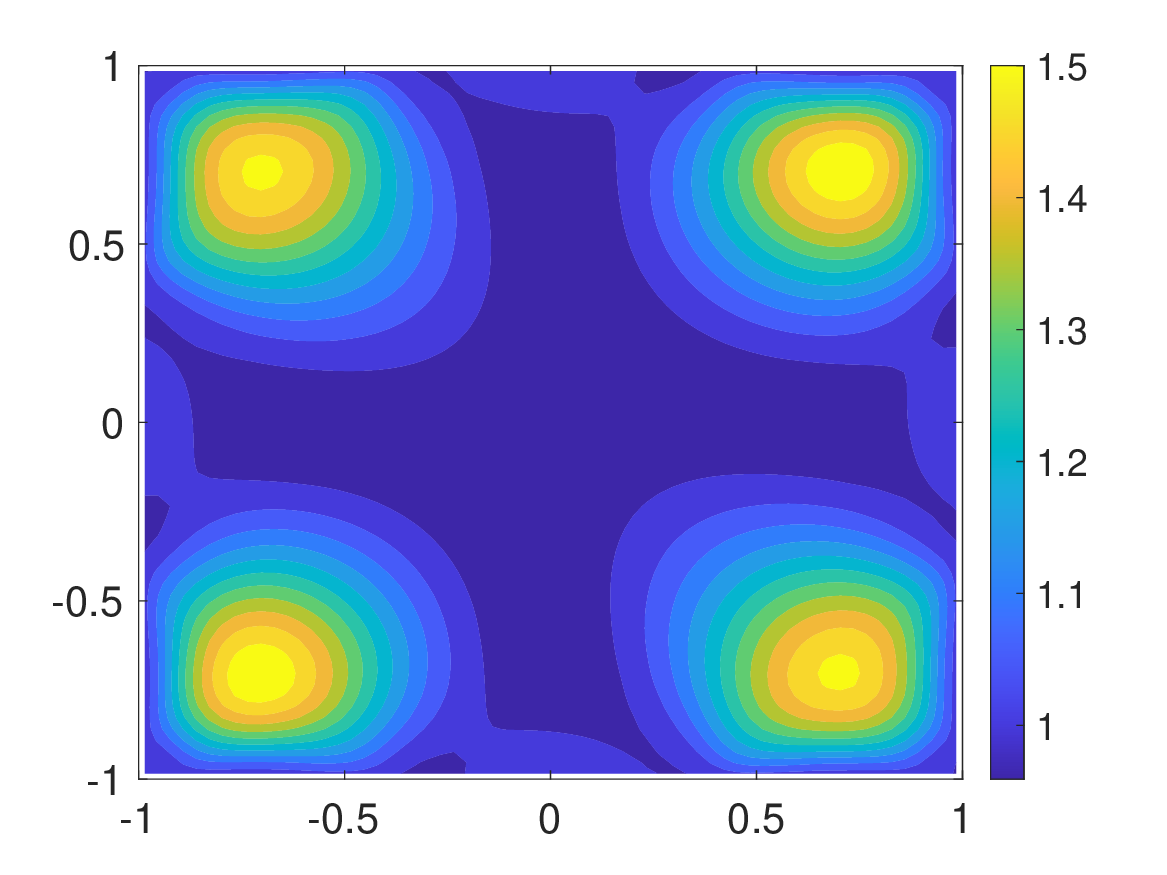}

\caption {Example \ref{ex:d2}:  computed conductivity with data noise
$\varepsilon=0.1\%$. Left: $16\times 16$; Middle: $32\times 32$; Right: $64\times 64$.}
\label{fig:e-3_7}
\end{center}
\end{figure}

{In this example, we consider a  discontinuous conductivity field with a constant background.
The true conductivity is shown in Figure \ref{fig:true7}. It has a height of 1.5 in the 
four squares and 1 anywhere else. A similar example can be found in \cite{WG}.
Figure \ref{fig:e-3_7}  shows the computed conductivity by MD-LDG with $P^2$ polynomial space of 
$16\times 16$, $32\times 32$ and $64\times 64$ from left to right with  data noise
$\varepsilon=0.1\%$.
We can see that the recoveries can well
capture the locations and heights of the four squares. 
We admit that due to the $H^1$ norm, the shapes of the conductivity have been smoothened somehow. }

\subsection{Example \ref{ex:d3}: Reconstruction of EIT: inclusions with a discontinuous background}
\label{ex:d3}
\begin{figure}
 \begin{center}
    \includegraphics[width = 2.in,height=1.7in] {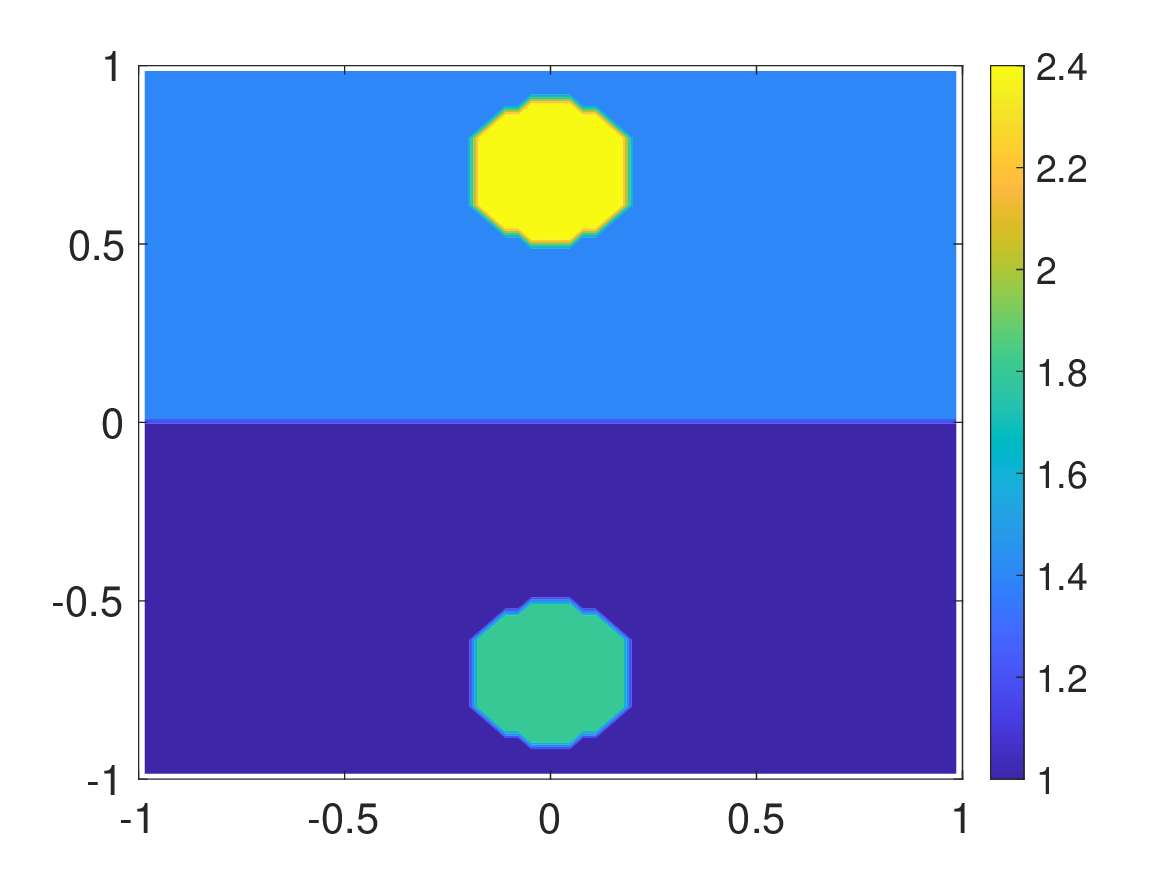}
\caption {Example \ref{ex:d3}:  true conductivity.}
\label{fig:true8}
\end{center}
\end{figure}

\begin{figure}
 \begin{center}
    \includegraphics[width = 2.in,height=1.7in] {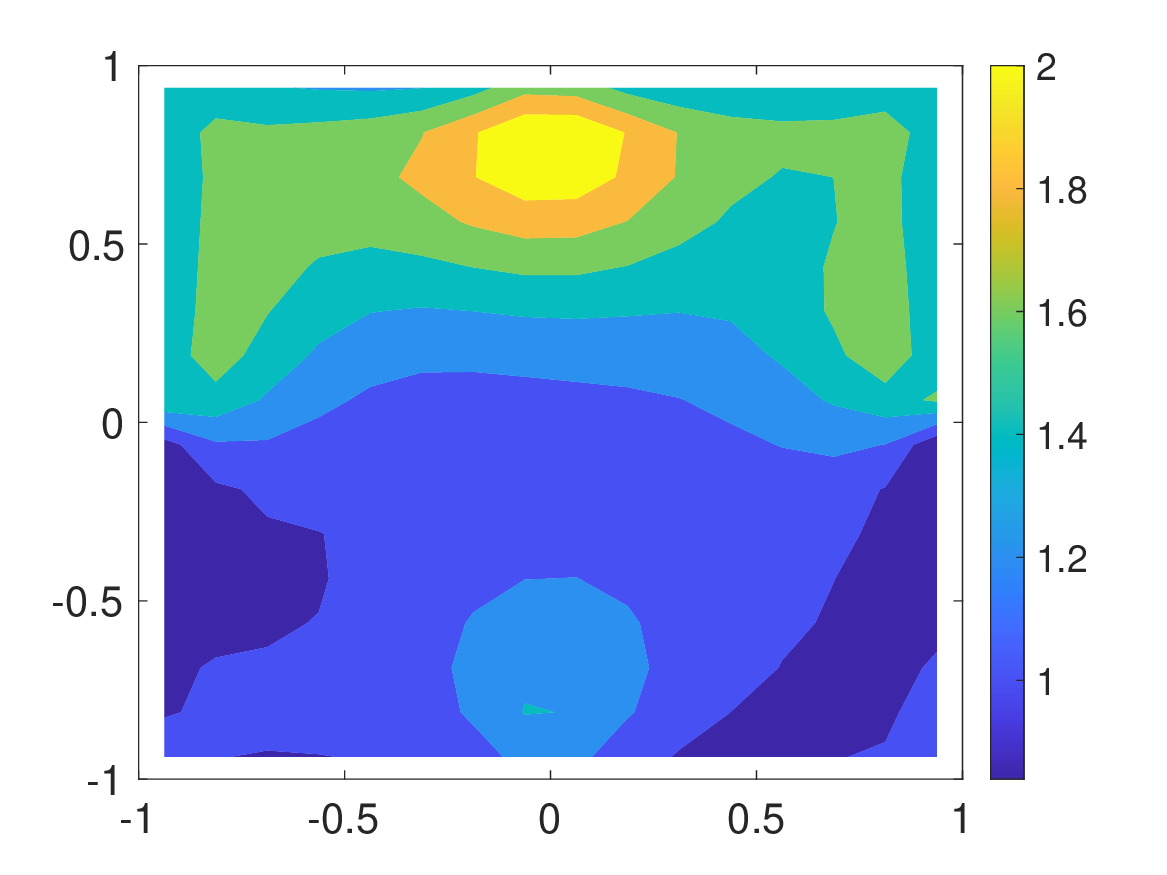}
        \includegraphics[width = 2.in,height=1.7in] {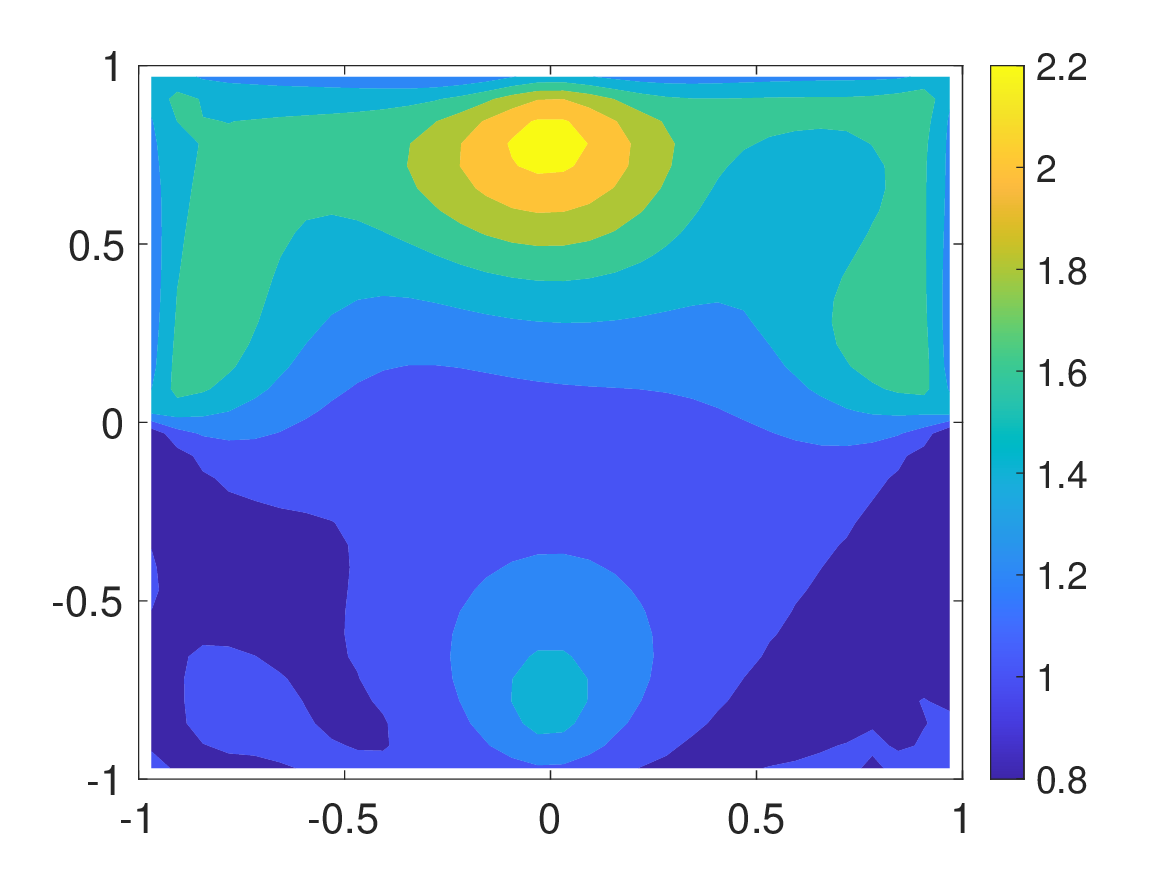}
   \includegraphics[width = 2.in,height=1.7in] {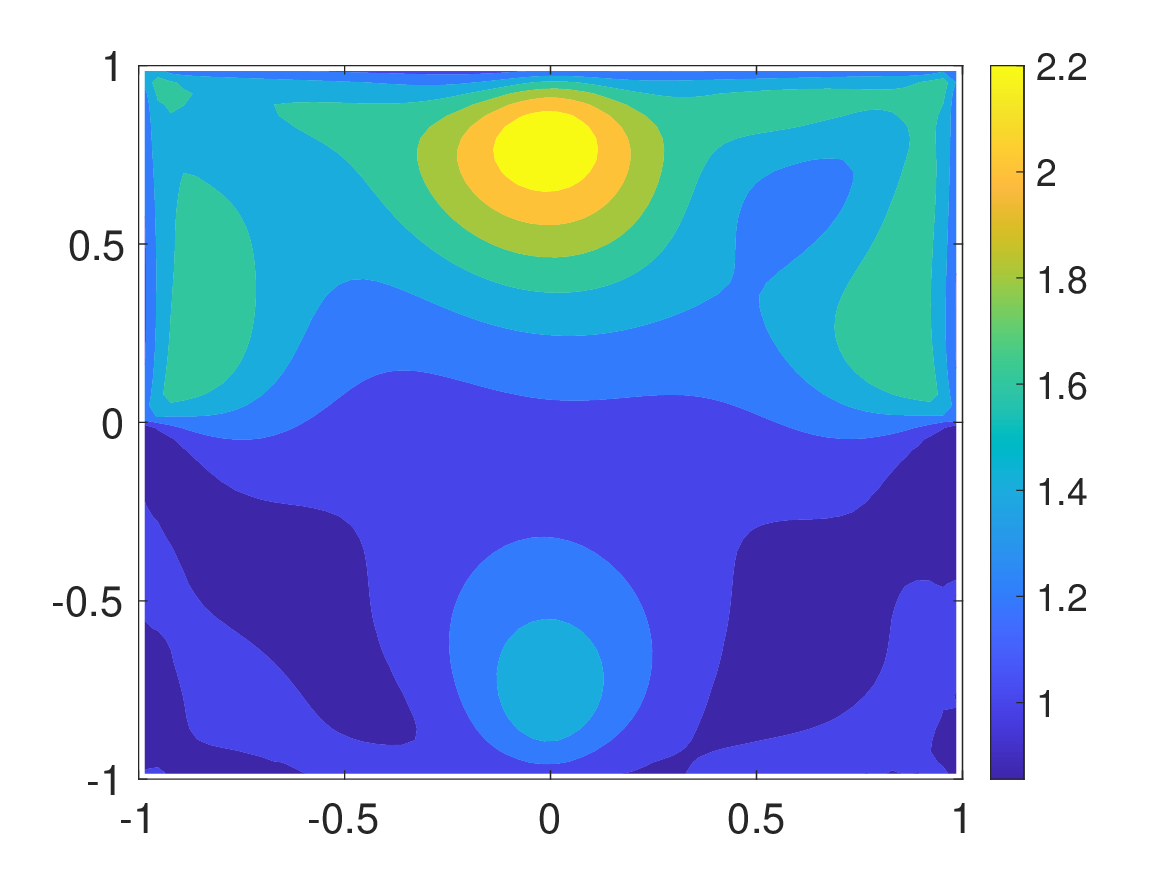}

\caption {Example \ref{ex:d3}:  computed conductivity with data noise
$\varepsilon=0.1\%$. Left: $16\times 16$; Middle: $32\times 32$; Right: $64\times 64$.}
\label{fig:e-3_8}
\end{center}
\end{figure}

{In the last example, we consider a  discontinuous conductivity field with a discontinuous background.
The background has a discontinuity at $y=0$ with a value of $1.5$ for $y>0$ and $1$ for $y<0$.
A similar example can be found in \cite{Jin12}.
The true conductivity consists of two circles centered at (0, 0.7) and (0, -0.7) with a height of 2.5 and 2, respectively, which is shown in Figure \ref{fig:true8}.
Figure \ref{fig:e-3_8}  shows the computed conductivity by MD-LDG with $P^2$ polynomial space of 
$16\times 16$, $32\times 32$ and $64\times 64$ from left to right with  data noise
$\varepsilon=0.1\%$.
We can see that  the recoveries can well
capture the locations, heights, and shapes of the two circles. }

\section{Concluding remarks}
In this paper, we consider the numerical reconstruction of the conductivity from Dirichlet-to-Neumann map. It is somehow different from the reconstruction from Neumann-to-Dirichlet map, where the latter has been extensively studied in the computational aspect. We developed a high order numerical method  for solving the Electrical Impedance Tomography problem which uses
a third order minimal-dissipation local discontinuous Galerkin method as the forward solver.
The reconstruction is based on the iterative least-squares method with Tikhonov regularization.
The efficiency and convergence of the algorithm are demonstrated by several numerical experiments \ww{including continuous and discontinuous conductivities}. The results show the proposed method can well recover the locations and shapes.

We remark that there are many details of the scheme that can be improved and investigated in future work.
In the present work, we consider the traditional $H^1$ penalty term in the regularization. We plan to work on other types of penalty terms including $l^1$ penalty for sparsity and total variation for discontinuity. 
We will also work on EIT problem with complete electrode model.
Finally, we will extend the proposed reconstruction method to other types of inverse problems.


\section*{Acknowledgements} 
We thank the reviewers for the valuable comments and suggestions.
This research did not receive any specific grant from funding agencies in the public, commercial, or
not-for-profit sectors.

\end{document}